\UseRawInputEncoding
\documentclass[12pt]{article}
\usepackage{booktabs,caption}
\usepackage[labelfont=bf]{caption}

\usepackage{changepage}
\usepackage{geometry}
\usepackage{amsfonts,color}
\definecolor{r}{rgb}{1.00,0.00,0.00}
\usepackage{amssymb,titlesec}
\titleformat*{\section}{\large\bf}
\titleformat*{\subsection}{\normalsize\bf}
\usepackage{lineno,hyperref}

\usepackage{amsfonts,amsmath,amsthm}
\usepackage[mathscr]{eucal}
\usepackage{graphicx}
\usepackage{epstopdf}
\usepackage{subfigure}
\usepackage{autobreak}
\allowdisplaybreaks

\newtheorem{thm}{Theorem}

\newtheorem{lem}[thm]{Lemma}
\newtheorem{cory}[thm]{Corollary}
\newtheorem{rmk}{Remark}[section]

\newcommand{\RNum}[1]{\uppercase\expandafter{\romannumeral #1\relax}}

\renewcommand{\baselinestretch}{1}

\begin{document}
\numberwithin{equation}{section}
\vskip 5mm

\noindent
\begin{Large}
 \centerline{\textbf{Spectrally negative L\'{e}vy risk model under mixed}}\\
 \centerline{\textbf{ratcheting-periodic dividend strategies }}
 \end{Large}
\vskip 3mm

\noindent
\begin{adjustwidth}{1cm}{0cm}\centerline{SUN FUYUN$^{a}$\ \ AND SONG ZHANJIE$^{a}$}\ \
\begin{small}
\vskip 3mm
\noindent \centerline{${}^{a}$School of Mathematics, Tianjin University, Tianjin 300350, China}\\
\end{small}
\end{adjustwidth}

\noindent
\begin{adjustwidth}{1cm}{0cm}
\begin{small}
\textbf{ABSTRACT} \hspace{2mm} In this paper, we consider the mixed ratcheting-periodic dividend strategies for spectrally negative L\'{e}vy risk model, in which dividend payments can both be made continuously without falling and discretely at the jump times of an independent Poisson process. The expected net present value(NPV) of dividends paid up to ruin and the Laplace transform of the ruin time are obtained by using L\'{e}vy fluctuation theory. All the results are expressed in terms of scale functions. Finally, numerical results for Brownian motion with drift are given.
\vskip 3mm
\noindent
\textbf{KEYWORDS} \hspace{2mm} Expected net present value(NPV) of dividends; Ratcheting dividend strategy; Periodic dividend strategy; Spectrally negative L\'{e}vy process; Scale function; Laplace transform.
\vskip 3mm
\noindent
\textbf{MATHEMATICS SUBJECT CLASSIFICATION} \hspace{2mm} 91B30, 97M30, 60J75.
\end{small}
\end{adjustwidth}
\vskip 4mm

\begin{small}
Address correspondence to Fuyun Sun, School of Mathematics, Tianjin University, Tianjin 300350, China; E-mail: sunfy@tju.edu.cn

\end{small}

\vskip 4mm
\noindent

\section{Introduction}\label{Introduction}

In actuarial risk theory, the L\'{e}vy risk models with barrier dividend strategy have been studied extensively, see, e.g. Loeffen (2008); Kyprianou and Loeffen (2010); Kyprianou, Loeffen, and P\'{e}rez (2012); Yin and Wen (2013); Yin, Shen, and Wen (2013); Yin, Wen, and Zhao (2014), among others.
In reality, the reduction in the dividend rate may have a negative psychological impact on shareholders, which may lead to a decrease in earnings.
In order to avoid the above situation, a ratcheting dividend strategy in risk theory was considered
by Albrecher, B\"{a}uerle, and Bladt (2018), where the dividend rate can never decrease. They obtained the
expected value of the aggregate discounted dividend payments until ruin under a ratcheting
strategy for a L\'{e}vy risk model. After that there are some recent papers on ratcheting
dividend strategy in different risk models Zhang and Liu (2020); Albrecher, Azcue, Muler (2020a,b); Song and Sun (2021).

However, for companies, the board should check the balance firstly on a periodic basis and then decide the dividend payments paid to the shareholders, turning out either continuous payment streams or one-off dividend payments at discrete time points. Hence the board of the company may choose a combined dividend strategy. Recently, mixed strategies have been studied by a lots of contributors in vary risk models, such as  Avanzi, Tu, and Wong (2020b) studied hybrid continuous and periodic barrier strategies in the dual model, Liu, Chen, and Hu (2020) considered threshold and periodic dividend strategies in a dual model with diffusion,  Avram, P\'{e}rez, and Yamazaki (2018) studied Parisian reflection below and classical reflection above in spectrally negative L\'{e}vy processes and other papers Zhang and Han (2017); Dong and Zhou (2019); Avanzi, Tu, and Wong (2016); Avanzi et al. (2017); P\'{e}rez
and Yamazaki (2018), and so on. 
Motivated by the board's reasonable behaves and these works, we consider in this paper the spectrally negative L\'{e}vy case with mixed ratcheting-periodic dividend strategy, which is a combination of continuous ratcheting dividend strategy and discrete periodic dividend strategy. The ratcheting dividend rate we considered during the lifetime can be increased once to a higher level from the original level, like in Albrecher, B\"{a}uerle, and Bladt (2018). The periodic dividend strategy we considered can be portrayed as that the surplus process is pushed down to a preset barrier whenever it is above the barrier at the periodic dividend decision times.

The expected net present value(NPV) of dividends paid up to ruin and the Laplace
transform of the ruin time have been extensively studied in the literature by
using the resolvent measure, the Laplace transform of the occupation times and
other fluctuation identities. We also obtained the numerical optimal barriers value under the special model, because the analysis solution of the optimal barriers is more complicated. For more studies on the expected NPV of dividends and
the Laplace transform of the ruin time,
see Yin and Yuen (2011); Shen, Yin, anf Yuen (2013); Dong, Yin, and Dai (2019); Avanzi, Lau, and Wong (2021); Avanzi, Tu, and Wong (2020a); Li et al. (2021), etc.
It is worth mentioning that, ratcheting-periodic dividend strategy can reduce to the pure periodic barrier strategy(see e.g. Avram, P\'{e}rez, and Yamazaki (2018)) and to the pure ratcheting dividend strategy(see e.g. Albrecher, B\"{a}uerle, and Bladt (2018)) under certain conditions.

The rest of the paper is organized as follows. In Section $\ref{Pre}$, we list and recall some preliminaries: spectrally negative L\'{e}vy risk processes in Subsection $\ref{Sper}$, some associated scale functions in Subsection $\ref{Review}$, L\'{e}vy risk models with periodic barrier strategy in Subsection $\ref{Lev}$, and the definition of the ratcheting-periodic dividend strategy and the construction of the corresponding controlled surplus process in Subsection $\ref{Mixed}$. The expressions of the expected NPV of dividends up to ruin and the Laplace transform of the ruin time are discussed in Sections $\ref{The}$ and $\ref{Lap}$, respectively. Section $\ref{Num}$ shows the analysis with Brownian motion. A conclusion is given in Section $\ref{Con}$.

\section{Preliminaries}\label{Pre}
\subsection{Spectrally negative L\'{e}vy processes}\label{Sper}
Let us consider a spectrally negative L\'{e}vy process $Y=\{Y(t),~t\ge 0\}$, i.e. a L\'{e}vy process with only negative jumps. We assume that $Y(0)= y$ and the drift of this process is positive.
For $y\in \mathbb{R}$, we denote by $\mathbb{P}_y$ the law of $Y$ when it starts at $y$. Accordingly, I shall write $\mathbb{E}_y$ for the associated expectation operator. Throughout this work define the Laplace exponent $\psi(\theta)=\log\mathbb{E}(e^{\theta Y(1)})$, which is finite for at least all $\theta\ge 0$, by the L\'{e}vy-Khintchine formula(see e.g.  Kuznetsov, Kyprianou,and Rivero (2013); Chan, Kyprianou, and Savov (2011); Kyprianou, Loeffen,and P\'{e}rez (2012)): 
\begin{align*}
\psi(\theta)=\frac{1}{2}\sigma^2\theta^2
+\gamma\theta+\int_{(-\infty, 0)}(e^{\theta y}-1-\theta y \textbf{1}_{\{y>-1\}})\Pi(dy)
,~~\theta\ge 0,
\end{align*}
where $\gamma\in\mathbb{R}$, $\sigma\ge 0$ and $\Pi$ is a measure on $(-\infty, 0)$ called the L\'{e}vy measure of $Y$ that satisfies
\begin{align*}
\int_{(-\infty, 0)}(1\wedge y^2)\Pi(dy)<\infty.
\end{align*}
It is well-known that $Y$ has paths of bounded variation if and only if $\sigma= 0$ and $\int_{(-1, 0)}|y|\Pi(dy)<\infty$; in this case, $Y$ can be written as
$Y(t) = ct -  S(t),~t \ge 0 $,
where
\begin{align}\label{2.1}
c = \gamma - \int_{(-1, 0)}y\Pi(dy)
\end{align}
and
$\{S(t), t \ge 0\}$ is a driftless  subordinator. Note that necessarily
$c > 0$, since we have ruled out the case that $Y$ has monotone paths; its Laplace exponent is given by
\begin{align*}
\psi(\theta)=c\theta+\int_{(-\infty, 0)}(e^{\theta y}-1)\Pi(dy), ~~~\theta\ge 0.
\end{align*}
Throughout the paper, we assume that $\mathbb{E}[Y(1)]=\psi^{'}(0+)<\infty.$

\subsection{Review of scale functions}\label{Review}

Let $X(t) = Y(t) - c_1 t$, and $\widetilde{X}(t)=Y(t) - c_1 t - c_2 t$,~$t\ge 0$, $c_1\ge 0$,~$c_2>0$. In this paper, We also assume that $c_1+c_2\leq c$(well defined in equation $(\ref{2.1})$) for ensuring both processes $X=\{X(t),~t\ge 0\}$ and $\widetilde{X}=\{\widetilde{X}(t),~t\ge 0\}$ have positive drift. To avoid confusion of the notations, we assume that $\psi(\theta)$ denotes the Laplace exponent of process $X$ for the rest of the paper.
Denote by $\Phi( \delta )$ the largest
root of the equation $\psi(\theta)=\delta$, $i.e.$, $\Phi( \delta )=\sup\{\theta\ge0 : \psi(\theta)-\delta=0\}$.
We now recall the definition of the $\delta$-scale function $W^{(\delta)}(x)$.
For each $\delta\ge 0$, there exists a continuous and increasing function $W^{(\delta)} : \mathbb{R}\rightarrow [0, \infty)$, which is called the $\delta$-scale function of the process $X$.
The $\delta$-scale function $W^{(\delta)}(x)$ and some related functions of the process $X$ are defined in such a way:
\begin{align*}
\int_0^{\infty}e^{-\theta x}W^{(\delta)}(x)dx=
\frac{1}{\psi(\theta)-\delta},~~\theta>\Phi(\delta),
\end{align*}
and for $\delta\ge 0$ and $x\in \mathbb{R}$,
\begin{align*}
&\overline{W}^{(\delta)}(x)=\int_0^{x}W^{(\delta)}(u)du,~~
\overline{\overline{W}}^{(\delta)}(x)=\int_0^{x}\int_0^{y}W^{(\delta)}(u)dudy,~\\
&Z^{(\delta)}(x)=1+\delta \overline{W}^{(\delta)}(x),~~~~~
\overline{Z}^{(\delta)}(x)=\int_0^{x}Z^{(\delta)}(u)du.
\end{align*}
Noting that $W^{(\delta)}(x)=0$ for $-\infty<x<0$, we then have
$\overline{W}^{(\delta)}(x)=0$, $\overline{\overline{W}}^{(\delta)}(x)=0$, $Z^{(\delta)}(x)=1$ and
$\overline{Z}^{(\delta)}(x)=x$, respectively, for $x \le 0$.

Define also
\begin{align*}
Z^{(\delta)}(x,\theta)=e^{\theta x}\left(1+(\delta-\psi(\theta))\int_0^{x}e^{-\theta u}W^{(\delta)}(u)du \right),~\theta\ge 0, ~x\in \mathbb{R},
\end{align*}
and its partial derivative with respect to $x$:
\begin{align*}
Z^{(\delta)'}(x,\theta)=\theta Z^{(\delta)}(x,\theta)+
(\delta-\psi(\theta))W^{(\delta)}(x),~\theta\ge 0, ~x\in \mathbb{R}.
\end{align*}
In particular, for $x\in \mathbb{R}$, $Z^{(\delta)}(x, 0)= Z^{(\delta)}(x)$ and for $\delta\ge 0$,
\begin{align*}
&Z^{(\delta)}(x,\Phi(\delta+\gamma))=
e^{ \Phi(\delta+\gamma) x}
\left(1-\gamma\int_0^{x}e^{-\Phi(\delta+\gamma) u}
W^{(\delta)}(u)du \right).
\end{align*}

We give some more notations, which will be used later:
for any measurable function $f : \mathbb{R}\rightarrow\mathbb{R}$,
\begin{align*}
\mathcal{M}_b^{(\delta, \gamma)}f(x):=f(x+b)+
\gamma\int_0^x W^{(\delta+\gamma)}(x-y)f(y+b)dy,~b>0,~x\in\mathbb{R}.
\end{align*}
In particular, we let, for $b>0$, $\delta \ge 0$ and $x\in\mathbb{R}$,
\begin{align*}
W_b^{(\delta, \gamma)}(x):=
\mathcal{M}_b^{(\delta, \gamma)}W^{(\delta)}(x),~~~ Z_b^{(\delta, \gamma)}(x):=
\mathcal{M}_b^{(\delta, \gamma)}Z^{(\delta)}(x),
\end{align*}
which are taken from P\'{e}rez and Yamazaki (2018). Note that similar generalized scale functions are also introduced in Avram, P\'{e}rez, and Yamaziki (2018).

Define also for $b>0$, $\delta \ge 0$ and $x\in\mathbb{R}$,
\begin{align*}
&I_b^{(\delta, \gamma)}(x):=\frac{W_b^{(\delta, \gamma)}(x)}{W^{(\delta)}(b)}
-\gamma\overline{W}^{(\delta+\gamma)}(x),\\
&J_b^{(\delta, \gamma)}(x):=Z_b^{(\delta, \gamma)}(x)-
\gamma Z^{(\delta)}(b)\overline{W}^{(\delta+\gamma)}(x).
\end{align*}
Note in particular that for $-b<x<0$,
\begin{align}\label{1}
I_b^{(\delta, \gamma)}(x)=\frac{W_b^{(\delta )}(x+b)}{W^{(\delta)}(b)},~~
J_b^{(\delta, \gamma)}(x)=Z_b^{(\delta )}(x+b).
\end{align}

The corresponding functions for the process $\widetilde{X}$ will be denoted by $\mathbb{W}^{(\delta)}(x)$, $\overline{\mathbb{W}}^{(\delta)}(x)$,
$\overline{\overline{\mathbb{W}}}^{(\delta)}(x)$, $\mathbb{Z}^{(\delta)}(x)$, $\overline{\mathbb{\mathbb{Z}}}^{(\delta)}(x)$, $\mathbb{Z}^{(\delta)}(x,\theta)$, $\mathbb{W}_b^{(\delta, \gamma)}(x)$,
$\mathbb{Z}_b^{(\delta, \gamma)}(x)$,
$\mathbb{I}_b^{(\delta, \gamma)}(x)$ and $\mathbb{J}_b^{(\delta, \gamma)}(x)$,
  respectively, and $\phi(\delta)$ shall be the corresponding largest root of equation $\psi(\theta)-c_2 \theta=\delta$.

\subsection{L\'{e}vy risk models with periodic barrier strategy}\label{Lev}
%

Firstly let $\{T_i\}_{i=1}^{\infty}\left(T_1<T_2<\cdots\right)$ be an increasing sequence of jump times of an independent Poisson process with rate $\gamma(\gamma>0)$. Whenever the observed surplus level at $T_i$ is larger than $a(a\ge 0)$, the excess value will be paid off as dividend. We construct the L\'{e}vy risk model with periodic barrier strategy $X^a=\{X^a(t),~t\ge 0\}$ as follows.
Specifically, we have:
\begin{align*}
 X^a(t)=X(t),~0\le t\le W^{+}_0(1),
\end{align*}
where:
\begin{align}\label{2.2}
W^{+}_0(1):=\inf \{T_i~:~X(T_i)>a\};
\end{align}
here and throughout, let $\inf \emptyset = \infty$. The process then jumps downward by
$X(W^{+}_0(1))-a$ so that $X^a(W^{+}_0(1))=a$. For $W^{+}_0(1)\le t< W^{+}_0(2):=\inf \{T_i>W^{+}_0(1):X^a(W^{+}_0(1)-)>a\}$, we have $X^a(t)=X(t)-X(W^{+}_0(1))$, and $X^a(W^{+}_0(2))=a$. The process can be constructed by repeating this procedure.

Suppose $D_P(t)$ is the cumulative amount of periodic dividends until time $t\ge 0$. Then, we have:
\begin{align}\label{2.3}
X^a(t)=X(t)-D_P(t),~t\ge 0,
\end{align}
with
$$D_P(t)=\sum_{W^{+}_0(i)\le t}\left(X^a( W^{+}_0(i)-)-a\right),~t\ge 0,$$
where $\{W^{+}_0(n),~n\ge 1\}$ can be constructed inductively by $(\ref{2.2})$ and:
$$W^{+}_0(n+1):=\inf \{T_i>W^{+}_0(n)~:~X^a(T_i-)>a\},~i\ge 1.$$

Similar to the construction method of process $X^a$, we have the process $\widetilde{X}^a=\{\widetilde{X}^a(t),~t\ge 0\}=\{\widetilde{X}(t)-\widetilde{D}_P(t),~t\ge 0\}$, where $\widetilde{D}_P(t)$ denotes the corresponding cumulative periodic dividends until time $t\ge 0$.

For the processes
$X^a $ and
$\widetilde{X}^a $, define for a fixed $a\ge0$ the first passage times:
\begin{align*}
&\tau^{-}_{0} :=\inf \{t\ge 0:X^a(t)<0\},~~
\tau^{+}_{a} :=\inf\{t\ge 0:X^a(t)>a\},\\
&\widetilde{\tau}^{-}_{0} :=\inf \{t\ge 0:\widetilde{X}^a(t)<0\},~~
\widetilde{\tau}^{+}_{a} :=\inf\{t\ge 0:\widetilde{X}^a(t)>a\},
\end{align*}
with the usual convention $\inf{\emptyset}=\infty$.

\subsection{Mixed ratcheting-periodic dividend strategies}\label{Mixed}

In this subsection, motivated by those works Albrecher, B\"{a}uerle, and Bladt (2018),  Liu, Chen, and Hu (2020) and P\'{e}rez and Yamazaki (2018),~
we propose a mixed dividend strategy, which contains a ratcheting dividend strategy and a periodic dividend strategy.
With the mixed dividend strategy, we modify the process $Y$ as follows. Whenever the observed surplus level at $T_i$(well defined in subsection \ref{Lev}) is larger than $a$, the excess value will be paid off as dividend. During the internal times $(T_{i-1},~T_i)$, dividends are paid at a fixed constant rate $c_1$ until the first time when the surplus process hits a predetermined barrier $b(b\ge a)$ and from this point the dividend rate will be ratcheted to $c_1+c_2$($c_1+c_2 \le c$) from a fixed constant $c_2$ and stays at this higher level until ruin. In order to give the mathematical descriptions of the modified surplus process
$\widetilde{Y}^{(a,b)}=\left\{\widetilde{Y}^{(a,b)}(t),t\geq 0\right\}$, starting from $\widetilde{Y}^{(a,b)}(0)=Y(0)=y$, under the mixed ratcheting and periodic dividend strategies, we define an auxiliary process $\widetilde{Y}=\left\{\widetilde{Y}(t),t\ge0\right\}$ as follows:
$$ d\widetilde{Y}(t)=dY(t)-c_1 dt - c_2 I_{\{M(t)>b\}}dt,~t\ge 0,$$
where $M(t)=\sup_{0\le s\le t}Y(s)$.
Then the modified process $\widetilde{Y}^(a,b)$ is given by:
$$ \widetilde{Y}^{(a,b)}(t)=\widetilde{Y}(t),~~0\le t < T^+_a(1),$$
where
\begin{align}\label{2}
T^+_a(1):=\inf\left\{T_i:\widetilde{Y}(T_i)>a\right\}.
\end{align}
The process then jumps downward by
$\widetilde{Y}(T^+_a(1))-a$
so that $\widetilde{Y}^{(a,b)}(T^+_a(1))=a$.
For
$T^+_a(1)\le t< T^+_a(2):=\inf\left\{T_i>T^+_a(1):\right.$
$\left.
\widetilde{Y}^{(a,b)}(T_i-)>a
\right\}$, we have
$\widetilde{Y}^{(a,b)}(t)=
\widetilde{Y}(t)-\left(\widetilde{Y}(T^+_a(1))-a\right)$, and $\widetilde{Y}^{(a,b)}(T^+_a(2))=a$.
The process
$\widetilde{Y}^{(a,b)}$
can be constructed by repeating this procedure.

Without loss of generality, we set $T_0=0$. Note that $T_0$ is not a dividend decision time, therefore $\widetilde{Y}^{(a,b)}(0)=y$ even if $y>a$. Then we have
\begin{align}\label{3}
\widetilde{Y}^{(a,b)}(t) =
\widetilde{Y}(t) - L_R(t)-L_P(t),~t \ge 0,
\end{align}
with
\begin{align*}
&L_R (t) :=\int_{[0,t]} e^{-\delta s}\left(c_1+c_2I_{\{M(s)>b\}}\right)ds,~~t \ge 0,\\
&L_P (t) :=\sum_{T^+_a(i) \le t}e^{-\delta T^+_a(i)}
\left(\widetilde{Y}^{(a,b)}(T^+_a(i)-)-a\right),~~t\ge 0,
\end{align*}
where $\left\{T^+_a( n ); n \ge 1\right\}$ can be constructed inductively by $(\ref{2})$ and
\begin{align*}
T^+_a(n + 1):= \inf \left\{ T_i > T^+_a(n)
: \widetilde{Y}^{(a,b)}(T_i-)>a \right\},~n \ge 1 .
\end{align*}

For the process
$\widetilde{Y}^{(a,b)}$, define for a fixed $b>0$ the first passage times:
\begin{align*}
\sigma^{-}_{0} :=\inf \{t\ge 0:\widetilde{Y}^{(a,b)}t)<0\},~~~
\sigma^{+}_{b} :=\inf\{t\ge 0:\widetilde{Y}^{(a,b)}(t)>b\},
\end{align*}
with the usual convention $\inf{\emptyset}=\infty$.

Analyzing these processes $X$, $\widetilde{X}$, $X^a$, $\widetilde{X}^a$, and $\widetilde{Y}^{(a,b)}$, we can obtain the following lemma.
\begin{lem}\label{lem3.1}
\begin{align}
&\forall~t\le \sigma^-_0\wedge\sigma^+_b,~\widetilde{Y}^{(a,b)}(t)=X^a(t)=X(t)-D_P(t),~a.s.  \label{3.1}\\
&\forall~t\le \sigma^+_a,~\widetilde{Y}^{(a,b)}(t)=X^a(t)=X(t),~a.s.  \label{3.2}\\
&\forall~t> \sigma^+_b,~\widetilde{Y}^{(a,b)}(t)=\widetilde{X}^a(t)=\widetilde{X}(t)-\widetilde{D}_P(t),~a.s. \label{3.3}
\end{align}
\end{lem}
\section{The expected net present value of dividends}\label{The}
In this section, we present the expression of the expected NPV of dividends up to ruin via scale functions. In our assumption, the inter-dividend-decision times $W_i=T_i-T_{i-1}$, $i\ge 1$ are i.i.d. and exponentially distributed with mean $1/\gamma$ $(\gamma>0)$.
Define the expected NPV of dividends paid up to ruin by
\begin{align*}
V(y;a,b):&=\mathbb{E}\left[L_R(\widetilde{\tau}^{(a,b)})+L_P(\widetilde{\tau}^{(a,b)})\big{|}
\widetilde{Y}^{(a,b)}(0)=y \right]\\
&=
\mathbb{E}_y\left[L_R (\widetilde{\tau}^{(a,b)})+L_P(\widetilde{\tau}^{(a,b)}) \right],~~y\ge 0,
\end{align*}
where $\widetilde{\tau}^{(a,b)}=\inf\left\{t\ge 0: \widetilde{Y}^{(a,b)}(t)<0\right\}$ is the ruin time of the process $\widetilde{Y}^{(a,b)}$ well defined in $(\ref{3})$, with $\inf \emptyset =\infty$. Note that $\widetilde{\tau}^{(a,b)}=\sigma^-_0$.

According to our understanding and experience of the barrier dividend strategy,
$V(y;a,b)$ is different in the intervals $[0,a)$, $[a,b)$ and $[b,\infty)$ as both $a$ and $b$ are the ratcheting-periodic barriers. Then for easy identification, we denote by $V(y;a,b)=V_U(y;a,b)$ for $0<b\le y $, $V(y;a,b)=V_M(y;a,b)$ for $a\le y < b $ and $V(y;a,b)=V_M(y;a,b)$ for $0\le y < a $,
which will be given in Theorems $\ref{thm3.1}$--$\ref{thm3.3}$.

Before stating Theorems $\ref{thm3.1}$--$\ref{thm3.3}$, we first present two technical lemmas, which will be used to the proofs later. By directly applying Corollaries $1.$ $(ii)$ and $3.$ $(i)$ in the works of P\'{e}rez and Yamazaki (2018), we obtain the first Lemma $\ref{lem3.2}$:
\begin{lem}\label{lem3.2}
For $\delta\ge 0$, $a>0$ and $y\ge a$, we have
\begin{align}\label{3.4}
&\mathbb{E}_y\left[e^{-\delta \widetilde{\tau}^{-}_{0}}~;~\widetilde{\tau}^{-}_{0}<\infty \right]=
\mathbb{\mathbb{J}}_a^{(\delta,\gamma)}(y-a)-\mathbb{I}_a^{(\delta,\gamma)}(y-a)
\frac{\delta \mathbb{W}^{(\delta)}(a) \mathbb{Z}^{(\delta)}(a,\phi(\delta+\gamma))}{\mathbb{Z}^{(\delta)'}(a,\phi(\delta+\gamma))},
\end{align}
and
 \begin{align}\label{3.5}
 \mathbb{E}_y\left[\int_{[0,{\widetilde{\tau}^{-}_{0}}]}e^{-\delta t}d\widetilde{D}_P(t)
    \right]
 &=\mathbb{E}_y\left[\sum_{T_i\le \widetilde{\tau}^{-}_{0}}e^{-\delta T_i}\left(\widetilde{X}^{a}(T_i-)-a\right)\right]\nonumber\\
 &=
 \gamma\left[\frac{\mathbb{I}_{a}^{(\delta,\gamma)}(y-a)
\mathbb{W}^{(\delta)}(a)}{\phi(\delta+\gamma)
\mathbb{Z}^{(\delta)'}(a,\phi(\delta+\gamma))}
-\overline{\overline{\mathbb{W}}}^{(\delta+\gamma)}(y-a)  \right].
\end{align}
\end{lem}
Using the Theorems $1$, $2$ and $3$ in P\'{e}rez and Yamazaki (2018), we obtain the following Lemma $\ref{lem3.3}$.
\begin{lem}\label{lem3.3}
{For $\delta\ge 0$ and $b\ge a>0$, we have}
 \begin{align}
 &\mathbb{E}_y\left[e^{-\delta \tau^{+}_{b}}~;~\tau^{+}_{b}<\tau^{-}_{0} \right]
 =\frac{ {I}_{a}^{(\delta,\gamma)}(y-a)}{ {I}_{a}^{(\delta,\gamma)}(b-a)},~0\le y<b,\label{3.6}\\
 &\mathbb{E}_y\left[e^{-\delta \tau^{-}_{0}}~;~\tau^{-}_{0}<\tau^{+}_{b} \right]
 = {J}_{a}^{(\delta,\gamma)}(y-a)-\frac{ {I}_{a}^{(\delta,\gamma)}(y-a)}{ {I}_{a}^{(\delta,\gamma)}(b-a)}
 {J}_{a}^{(\delta,\gamma)}(b-a),~0\le y<b,\label{3.7}\\
&\mathbb{E}_y\left[\int_{[0,\tau^{+}_{b} \wedge \tau^{-}_{0}]}e^{-\delta t}d{D}_P(t)
    \right]
=\gamma \left(\overline{\overline{W}}^{(\delta+\gamma)}(b-a)
\frac{ {I}_{a}^{(\delta,\gamma)}(y-a)}{ {I}_{a}^{(\delta,\gamma)}(b-a)} -\overline{\overline{ W}}^{(\delta+\gamma)}(y-a) \right),~y\ge a.\label{3.8}
\end{align}
\end{lem}

Now we state the Theorems $\ref{thm3.1}$, $\ref{thm3.2}$, $\ref{thm3.3}$ and detail proofs as follows.
\begin{thm}\label{thm3.1} For $0<b\le y$, the expected NPV of
dividends paid up to ruin is given by
\begin{align}\label{8}
V_U(y;a,b)
= & \frac{c_1+c_2}{\delta}\left(1-\mathbb{J}_{a}^{(\delta,\gamma)}(y-a)
+\mathbb{I}_{a}^{(\delta,\gamma)}(y-a)
\frac{\delta\mathbb{W}^{(\delta)}(a)\mathbb{Z}^{(\delta)}(a,\phi(\delta+\gamma))}{\mathbb{Z}^{(\delta)'}(a,\phi(\delta+\gamma)) }\right)\nonumber\\
& + \gamma\left[\frac{\mathbb{I}_{a}^{(\delta,\gamma)}(y-a)
\mathbb{W}^{(\delta)}(a)}{\phi(\delta+\gamma)
\mathbb{Z}^{(\delta)'}(a,\phi(\delta+\gamma))}
-\overline{\overline{\mathbb{W}}}^{(\delta+\gamma)}(y-a)  \right].
\end{align}
\end{thm}

\textbf{Proof.} When the initial value $y \ge b$, dividends caused by ratcheting strategy are paid at rate $c_1+c_2$ until ruin. That means the whole modified risk model $\widetilde{Y}^{(a,b)}$ under mixed dividend strategy can be considered as the process $\widetilde{X}^a$ (see also equation (\ref{3.3}) in Lemma $\ref{lem3.1}$). But the total dividend amount of risk model $\widetilde{Y}^{(a,b)}$ includes the dividend amount generated by ratcheting strategy. Then we have
\begin{align}
 V_U(y;a,b)
 & =\mathbb{E}_y\left[\int_{[0,\widetilde{\tau}^{(a,b)}]} e^{-\delta t}
 \left(c_1+c_2 \right)dt+\sum_{T_a^+(i)\le \widetilde{\tau}^{(a,b)}}e^{-\delta T_a^+(i)}
 \left(\widetilde{Y}^{(a,b)}(T_a^+(i)-)-a\right)\right] \nonumber\\
 & =\mathbb{E}_y\left[\int_{[0,{\widetilde{\tau}_0^{-}}]} e^{-\delta t}
 \left(c_1+c_2 \right)dt\right]
 +\mathbb{E}_y\left[\int_{[0,{\widetilde{\tau}_0^{-}}]} e^{-\delta t}d\widetilde{D}_P(t)\right]\nonumber\\
& = \frac{c_1+c_2}{\delta}\left(1-\mathbb{E}_y\left[e^{-\delta \widetilde{\tau}^{-}_{0}}\right]\right)
+\mathbb{E}_y\left[\int_{[0,{\widetilde{\tau}_0^{-}}]} e^{-\delta t}d\widetilde{D}_P(t)\right]\label{9}.
\end{align}

Substituting equations $(\ref{3.4})$ and $(\ref{3.5})$ in Lemma $\ref{lem3.2}$ into $(\ref{9})$, we have $(\ref{8})$. This ends the proof.~ $\Box$

\begin{thm}\label{thm3.2} For $0<a\leq y<b$, the expected NPV of
dividends paid up to ruin is given by
\begin{align}\label{10}
V_M(y;a,b)
=& \frac{I_a^{(\delta,\gamma)}(y-a)}{I_a^{(\delta,\gamma)}(b-a)}
\left(V_U(b;a,b)+\frac{c_1}{\delta}\left(J_a^{(\delta,\gamma)}(b-a)-1 \right) +\gamma \overline{\overline{W}}^{(\delta+\gamma)}(b-a)\right)
\nonumber\\
&-\frac{c_1}{\delta}\left(J_a^{(\delta,\gamma)}(y-a)-1 \right)
-\gamma \overline{\overline{W}}^{(\delta+\gamma)}(y-a),
\end{align}
where
\begin{align*}
V_U(b;a,b)
= & \frac{c_1+c_2}{\delta}\left(1-\mathbb{J}_{a}^{(\delta,\gamma)}(b-a)
+\mathbb{I}_{a}^{(\delta,\gamma)}(b-a)
\frac{\delta\mathbb{W}^{(\delta)}(a)\mathbb{Z}^{(\delta)}(a,\phi(\delta+\gamma))}{\mathbb{Z}^{(\delta)'}(a,\phi(\delta+\gamma)) }\right)\nonumber\\
& + \gamma\left[\frac{\mathbb{I}_{a}^{(\delta,\gamma)}(b-a)
\mathbb{W}^{(\delta)}(a)}{\phi(\delta+\gamma)
\mathbb{Z}^{(\delta)'}(a,\phi(\delta+\gamma))}
-\overline{\overline{\mathbb{W}}}^{(\delta+\gamma)}(b-a)  \right]
\end{align*}
\end{thm}

\textbf{Proof}~~
Consider this case $a\leq y<b$, we discuss it in two ways. When the modified process $\widetilde{Y}^{(a,b)}$ does not reach $b$ before ruin, we can take it as process $X^a$ (see also equation $(\ref{3.1})$ in Lemma $\ref{lem3.1}$), where the dividends generated by the ratcheting strategy until ruin
$L_R (\widetilde{\tau}^{(a,b)})=c_1\int_{[0,\tau_0^-]} e^{-\delta t}dt$.
On the other hand, the modified process $\widetilde{Y}^{(a,b)}$ reaches $b$ before ruin.
In this case, the surplus $\widetilde{Y}^{(a,b)}$ can be described by equations $(\ref{3.1})$ and $(\ref{3.3})$. Then for this case we apply the strong Markov property at that point (up through $b$) in time, from which the process $X^a$ starting at $y$ dynamics change to the process $\widetilde{X}^a$ starting at $b$, which means ratcheting dividend rate changes to $c_1 + c_2$ from $c_1$.
We thus have
\begin{align}\label{11}
    V_M(y;a,b)
    & = c_1 \mathbb{E}_y\left[\int_{[0,\sigma_{0}^{-}\wedge \sigma_{b}^{+} ]}e^{-\delta t}dt\right]+
    \mathbb{E}_y\left[\int_{[0,\sigma_{0}^{-}\wedge \sigma_{b}^{+} ]}e^{-\delta t}dL_P(t)\right]\nonumber\\
    &~~ +\mathbb{E}_y\left[e^{-\delta\tau^{+}_{b}}~;~
    \tau^{+}_{b}<\tau^{-}_{0}\right]
    V_U(b;a,b).
\end{align}
By virtue of $(\ref{3.1})$ in Lemma $\ref{lem3.1}$, we get the first two terms of equation $(\ref{11})$
\begin{align}\label{12}
&  c_1 \mathbb{E}_y\left[\int_{[0,\sigma_{0}^{-}\wedge \sigma_{b}^{+} ]}e^{-\delta t}dt\right]+
    \mathbb{E}_y\left[\int_{[0,\sigma_{0}^{-}\wedge \sigma_{b}^{+} ]}e^{-\delta t}dL_P(t)\right]\nonumber\\
& = c_1 \mathbb{E}_y\left[\int_{[0,\tau_{0}^{-}\wedge \tau_{b}^{+} ]}e^{-\delta t}dt\right]+
    \mathbb{E}_y\left[\int_{[0,\tau_{0}^{-}\wedge \tau_{b}^{+} ]}e^{-\delta t}dD_P(t)\right]\nonumber\\
& =\frac{c_1}{\delta}\left(1-
\mathbb{E}_y\left[e^{-\delta\tau^{+}_{b}}~;~\tau^{+}_{b} <\tau^{-}_{0}\right]
-\mathbb{E}_y\left[e^{-\delta\tau^{-}_{0}}~;~\tau^{-}_{0} <\tau^{+}_{b} \right]\right)\nonumber\\
&~~~+\mathbb{E}_y\left[\int_{[0,\tau_{0}^{-}\wedge \tau_{b}^{+} ]}e^{-\delta t}dD_P(t)\right]\nonumber\\
&= \frac{c_1}{\delta}\left(1- {J}_{a}^{(\delta,\gamma)}(y-a)
+\left( {J}_{a}^{(\delta,\gamma)}(b-a)-1\right)
\frac{I_a^{(\delta,\gamma)}(y-a)}{I_a^{(\delta,\gamma)}(b-a)}
     \right)\nonumber\\
&~~~+\gamma\left(\overline{\overline{W}}^{(\delta+\gamma)}(b-a)
\frac{I_a^{(\delta,\gamma)}(y-a)}{I_a^{(\delta,\gamma)}(b-a)}
-\overline{\overline{W}}^{(\delta+\gamma)}(y-a)   \right)
.
\end{align}
The last equation is derived from $(\ref{3.6})$-$(\ref{3.8})$ in Lemma $\ref{lem3.3}$.
Applying equations $(\ref{12})$ and Lemma $\ref{lem3.2}$ to equation $(\ref{11})$, we have $(\ref{10})$.

Let $y=b$ in equation $(\ref{8})$, we then have the constant $V_U(b;a,b)$.
The proof is end. $\Box $%
%

\begin{thm}\label{thm3.3} For $0\leq y<a$, the expected NPV of
dividends paid up to ruin is given by
\begin{align}\label{a3}
V_L(y;a,b)
=& \frac{W^{(\delta)}(y)}{W^{(\delta)}(a)}
\left(V_M(a;a,b)+\frac{c_1}{\delta}(Z^{(\delta)}(a)-1)\right)
-\frac{c_1}{\delta}(Z^{(\delta)}(y)-1),
\end{align}
where
\begin{align*}
V_M(a;a,b)
=& \frac{1}{I_a^{(\delta,\gamma)}(b-a)}
\left(V_U(b;a,b)+\frac{c_1}{\delta}\left(J_a^{(\delta,\gamma)}(b-a)-1 \right) +\gamma \overline{\overline{W}}^{(\delta+\gamma)}(b-a)\right)
\nonumber\\&
-\frac{c_1}{\delta}\left(Z^{(\delta)}(a)-1 \right).
\end{align*}
\end{thm}
\textbf{Proof}~~ Consider this case $0\leq y<a$, we discuss it in two ways. When the modified process $\widetilde{Y}^{(a,b)}$ does not reach $a$ before ruin, we can take it as process $X$ (see also equation $(\ref{3.2})$ in Lamma $\ref{lem3.1}$), where the dividends generated by the ratcheting strategy until ruin $L_R (\widetilde{\tau}^{(a,b)})=c_1\int_{[0,\tau_0^-]} e^{-\delta t}dt$.
On the other hand, the modified process $\widetilde{Y}^{(a,b)}$ reaches $a$ before ruin.
Similar to the case of $0<a\leq y<b$, in this case we also apply the strong Markov property at that point (up through $a$) in time, from which the process $X$ or $X^a$ starting at $y$ dynamics change to the process $\widetilde{Y}^{(a,b)}$ with mixed ratcheting-periodic dividend strategy starting at $a$, which is the situation we discussed in the previous theorem.
We thus have
\begin{align}\label{a4}
    V_L(y;a,b)
    & = c_1 \mathbb{E}_y\left[\int_{[0,\sigma_{0}^{-}\wedge \sigma_{a}^{+} ]}e^{-\delta s}ds\right] +\mathbb{E}_y\left[e^{-\delta\tau^{+}_{a}}~;~
    \tau^{+}_{a}<\tau^{-}_{0}\right]
    V_M(a;a,b)\nonumber\\
    & = c_1 \mathbb{E}_y\left[\int_{[0,\tau_{0}^{-}\wedge \tau_{a}^{+} ]}e^{-\delta s}ds\right]+\mathbb{E}_y\left[e^{-\delta\tau^{+}_{a}}~;~
    \tau^{+}_{a}<\tau^{-}_{0}\right]
    V_M(a;a,b)\nonumber\\
    & = \frac{c_1}{\delta}
    \left(1- \mathbb{E}_y\left[e^{-\delta\tau^{+}_{a}}~;~\tau^{+}_{a} <\tau^{-}_{0}\right]
-\mathbb{E}_y\left[e^{-\delta\tau^{-}_{0}}~;~\tau^{-}_{0} <\tau^{+}_{a} \right]\right)\nonumber\\
&~~~ +\mathbb{E}_y\left[e^{-\delta\tau^{+}_{a} }~;~
    \tau^{+}_{a}<\tau^{-}_{0} \right]
    V_M(a;a,b).
\end{align}
By virtue of $(\ref{3.6})$ and $(\ref{3.7})$ in Lemma $\ref{lem3.3}$, we get the equation $(\ref{a3})$.

From equations $(\ref{1})$, we have $J_{a}^{(\delta,\gamma)}(0)={Z}^{(\delta)}(a)$ and
${I}_{a}^{(\delta,\gamma)}(0)$=1. Letting $y=a$ in equation $(\ref{10})$ and substituting these equations, we obtain the constant $V_M(a;a,b)$. The proof is end. $\Box $


Let $V^{a}(y;a),~y\ge 0,$ denote the expected NPV of dividends only under the periodic barrier strategy for spectrally negative L\'{e}vy  processes. By taking $c_1=0$ and $c_2~\downarrow~0$ in Theorems $\ref{thm3.1}$-$\ref{thm3.3}$, we have the following Corollary.

\begin{cory}\label{cory.1}~~For $\delta\ge 0$, $y\ge 0$ and $a>0$, we have
\begin{align*}
V^{a}(y;a)=\gamma\left[\frac{I_{a}^{(\delta,\gamma)}(y-a)
W^{(\delta)}(a)}{\phi(\delta+\gamma)
Z^{(\delta)'}(a,\phi(\delta+\gamma))}
-\overline{\overline{W}}^{(\delta+\gamma)}(y-a)  \right].
\end{align*}
\end{cory}
Corollary $\ref{cory.1}$ was used in Noba et al. (2018) to show the optimality of a periodic barrier strategy under the assumption that the L\'{e}vy measure has a completely monotone
density.

Let $V^R(y;b),~y\ge 0$ denote the expected NPV of dividends only under the ratcheting dividend strategy(single-rise) for spectrally negative L\'{e}vy  processes.
By taking $\gamma\downarrow 0$ in Theorems $\ref{thm3.1}$-$\ref{thm3.3}$, we have the following result, which coincides with the result in Theorem $2$ of Albrecher, B\"{a}uerle, and Bladt (2018).
\begin{cory}\label{cory.2}~~For $\delta, c_1, y\ge 0$ and $c_2, b>0$, we have
\begin{align*}
V^R(y;b)=
\left\{
  \begin{aligned}
  &\frac{c_1+c_2}{\delta}\left(1-\mathbb{Z}^{(\delta)}(y)+
\frac{\delta}{ \phi(\delta) }\mathbb{W}^{(\delta)}(y)\right),~~~~~~~~~0<b\le y,  \\
     & \frac{c_1+c_2}{\delta}\left(1-\mathbb{Z}^{(\delta)}(b)+
\frac{\delta}{ \phi(\delta) }\mathbb{W}^{(\delta)}(b)
\right)\frac{W^{(\delta)}(y)}{W^{(\delta)}(b)}\\
~&+ \frac{c_1}{\delta}\left(1-Z^{(\delta)}(y)+
(Z^{(\delta)}(b)-1) \frac{W^{(\delta)}(y)}{W^{(\delta)}(b)}\right),~~0\le y< b.
  \end{aligned}
  \right.
\end{align*}
\end{cory}

\textbf{Proof}~~It is obvious that $V^R(y;b)=\lim_{\gamma\downarrow 0}V^{R}(y;a,b)$. Next we aim to calculate this limit.

Note that for $\delta\ge 0$, $a>0$ and $y\in\mathbb{R}$,
\begin{align*}
  \lim_{\gamma\downarrow 0}\mathbb{I}_{a}^{(\delta,\gamma)}(y) =
  \frac{\mathbb{W}^{(\delta)}(y+a)}{\mathbb{W}^{(\delta)}(a)}~\rm{and}~
  \lim_{\gamma\downarrow 0}\mathbb{J}_{a}^{(\delta,\gamma)}(y)=
  \mathbb{Z}^{(\delta)}(y+a).
\end{align*}

According to the definition of $Z^{(\delta)}(a,\theta)$, we have
\begin{align*}
\lim_{\gamma\downarrow 0}\mathbb{Z}^{(\delta)}(a,\phi(\delta+\gamma))=
\mathbb{Z}^{(\delta)}(a,\phi(\delta))\end{align*}
as
$\lim_{\gamma\downarrow 0}\phi(\delta+\gamma)= \phi(\delta)$.
 Then
\begin{align*}
 &\lim_{\gamma\downarrow 0}\mathbb{Z}^{(\delta)'}(a,\phi(\delta+\gamma))\\
 &= \lim_{\gamma\downarrow 0}
 \left\{\phi(\delta+\gamma)\mathbb{Z}^{(\delta)}(a,\phi(\delta+\gamma))
 +[\delta-\psi(\phi(\delta+\gamma))+c_2 \phi(\delta+\gamma) ]W^{(\delta)}(a)  \right\} \\
 &= \phi(\delta)Z^{(\delta)}(a,\phi(\delta))
\end{align*}
as $\lim_{\gamma\downarrow 0}[\delta-\psi( \phi(\delta+\gamma) )+c_2 \phi(\delta+\gamma)]=\delta-\psi(\phi(\delta ))+c_2 \phi(\delta)=0  $.

From above equations and equations $(\ref{8})$, $(\ref{10})$, $(\ref{a3})$, we have
\begin{align*}
\lim_{\gamma\downarrow 0}V_U(y;a,b)&=\frac{c_1+c_2}{\delta}\left(1-\mathbb{Z}^{(\delta)}(y)+
  \frac{\mathbb{W}^{(\delta)}(y)}{\mathbb{W}^{(\delta)}(a)}\cdot
\mathbb{Z}^{(\delta)}(a,\phi(\delta)) \cdot \frac{\delta \mathbb{W}^{(\delta)}(a)}{\phi(\delta)\mathbb{Z}^{(\delta)}(a,\phi(\delta))}\right)\\
&~~+\lim_{\gamma\downarrow 0}\gamma \left(\frac{\mathbb{W}^{(\delta)}(y)}{\mathbb{W}^{(\delta)}(a)}\cdot
\frac{\mathbb{W}^{(\delta)}(a)}{\phi^2(\delta)\mathbb{Z}^{(\delta)}(a ,\phi(\delta))}
 - \overline{\overline{\mathbb{W}}}^{(\delta+\gamma)}(y-a) \right)\\
 &= \frac{c_1+c_2}{\delta}\left(1-\mathbb{Z}^{(\delta)}(y)+
  \frac{\delta}{\phi(\delta)}\mathbb{W}^{(\delta)}(y)\right),
\end{align*}
\begin{align*}
\lim_{\gamma\downarrow 0}V_M(y;a,b)
&=\lim_{\gamma\downarrow 0}
 \frac{I_a^{(\delta,\gamma)}(y-a)}{I_a^{(\delta,\gamma)}(b-a)}
\left(V_U(b;a,b)+\frac{c_1}{\delta}\left(J_a^{(\delta,\gamma)}(b-a)-1 \right) +\gamma \overline{\overline{W}}^{(\delta+\gamma)}(b-a)\right)
\nonumber\\
&~~-\lim_{\gamma\downarrow 0}\frac{c_1}{\delta}\left(J_a^{(\delta,\gamma)}(y-a)-1 \right)
-\lim_{\gamma\downarrow 0}\gamma \overline{\overline{W}}^{(\delta+\gamma)}(y-a)\\
&= \frac{W^{(\delta)}(y)}{W^{(\delta)}(b)}
\left(\lim_{\gamma\downarrow 0}V_U(b;a,b)+\frac{c_1}{\delta}\left(Z^{(\delta)}(b)-1 \right)\right)
-\frac{c_1}{\delta}\left(Z^{(\delta)}(y)-1 \right)\\
&=\frac{c_1+c_2}{\delta}\left(1-\mathbb{Z}^{(\delta)}(b)+
\frac{\delta}{\phi(\delta)}\mathbb{W}^{(\delta)}(b)
\right)\frac{W^{(\delta)}(y)}{W^{(\delta)}(b)}\\
&~~+ \frac{c_1}{\delta}\left(1-Z^{(\delta)}(y)+
(Z^{(\delta)}(b)-1) \frac{W^{(\delta)}(y)}{W^{(\delta)}(b)}\right),
\end{align*}
and
\begin{align*}
\lim_{\gamma\downarrow 0}V_L(y;a,b)
&=\lim_{\gamma\downarrow 0}
 \frac{W^{(\delta)}(y)}{W^{(\delta)}(a)}
\left(V_M(a;a,b)+\frac{c_1}{\delta}(Z^{(\delta)}(a)-1)\right)
-\frac{c_1}{\delta}(Z^{(\delta)}(y)-1)
\nonumber\\
&= \frac{W^{(\delta)}(y)}{W^{(\delta)}(a)}
\left\{\frac{c_1+c_2}{\delta}\left(1-\mathbb{Z}^{(\delta)}(b)+
\frac{\delta}{\phi(\delta)}\mathbb{W}^{(\delta)}(b)
\right)\frac{W^{(\delta)}(a)}{W^{(\delta)}(b)}\right.\\
&~~\left.+ \frac{c_1}{\delta}\left(1-Z^{(\delta)}(a)+
(Z^{(\delta)}(b)-1) \frac{W^{(\delta)}(a)}{W^{(\delta)}(b)}\right)\right\}\\
&~~+\frac{c_1}{\delta}(Z^{(\delta)}(a)-1)\frac{W^{(\delta)}(y)}{W^{(\delta)}(a)}
-\frac{c_1}{\delta}(Z^{(\delta)}(y)-1)
\nonumber\\
&=\frac{c_1+c_2}{\delta}\left(1-\mathbb{Z}^{(\delta)}(b)+
\frac{\delta}{\phi(\delta)}\mathbb{W}^{(\delta)}(b)
\right)\frac{W^{(\delta)}(y)}{W^{(\delta)}(b)}\\
&~~+ \frac{c_1}{\delta}\left(1-Z^{(\delta)}(y)+
(Z^{(\delta)}(b)-1) \frac{W^{(\delta)}(y)}{W^{(\delta)}(b)}\right).
\end{align*}
This ends the proof. $\Box$

\begin{rmk}\label{rmk.1}~~{\rm
 By taking $y \uparrow b$ in equation $(\ref{10})$ and taking $y \uparrow a$ in equation $(\ref{a3})$, we obtain
\begin{align*}
&\lim_{y \uparrow b }V_M(y;a,b)=V_U(b;a,b),\\
&\lim_{y \uparrow a }V_L(y;a,b)=V_M(a;a,b).
\end{align*}
This shows that for any spectrally negative L\'{e}vy  processes $Y(t)$, with respect to the initial value $y$, $V(y;a,b)\in C^0(0,\infty)$, even if $Y$ has paths of unbounded variation. The continuity property may be a necessary condition for existing of the optimal mixed strategy, if the mixed optimal strategy exists under some certain assumptions.
}
\end{rmk}

\section{ Laplace transform of the ruin time}\label{Lap}
In this section, we consider the distribution of the ruin time $\widetilde{\tau}^{(a,b)}$, under a ratcheting-periodic strategy.
For a fixed $\delta\ge 0$, define
\begin{align*}
L(y;a,b)=\mathbb{E}_y[e^{-\delta\widetilde{\tau}^{(a,b)}}~;~\widetilde{\tau}^{(a,b)}<\infty],  ~~b\ge a>0,~y\ge 0,
\end{align*}
as the expected NPV of a payment of $1$ at
the ruin time, and the Laplace transform of the probability density function of $\widetilde{\tau}^{(a,b)}$. We discuss how to calculate $L(y;a,b)$ as follows
for three cases $0\le y<a$, $a\le y<b$ and $y\ge b$, respectively.
\begin{thm}\label{thm4.1} For $\delta, c_1, y\ge 0$, $c_2>0$, and $b\ge a>0$, we have 
\begin{align*}
L(y;a,b)=
\left\{
  \begin{aligned}
  &\mathbb{J}_{a}^{(\delta,\gamma)}(y-a)
-\mathbb{I}_{a}^{(\delta,\gamma)}(y-a)
\frac{\delta\mathbb{W}^{(\delta)}(a)
\mathbb{Z}^{(\delta)}(a,\phi(\delta+\gamma))}{\mathbb{Z}^{(\delta)'}(a,\phi(\delta+\gamma))},
~~0<b\le y,  \\
&Z^{(\delta)}(y)+\frac{W^{(\delta)}(y)}{W^{(\delta)}(b)}
\left\{\mathbb{J}_{a}^{(\delta,\gamma)}(b-a)
-Z^{(\delta)}(b)\right.\\&
\left.~~~~~~~~~~~~-\mathbb{I}_{a}^{(\delta,\gamma)}(b-a)
\frac{\delta\mathbb{W}^{(\delta)}(a)\mathbb{Z}^{(\delta)}(a,\phi(\delta+\gamma))}{\mathbb{Z}^{(\delta)'}(a,\phi(\delta+\gamma))}
\right\},~0\le y<b.
  \end{aligned}
  \right.
\end{align*}
\end{thm}

\textbf{Proof}~~Consider the first case $y\ge b$, according to the Lemma $\ref{lem3.1}$, it is easy to obtain that
\begin{align}\label{14}
L(y;a,b)&=\mathbb{E}_y[e^{-\delta\widetilde{\tau}^{(a,b)}}~;~\widetilde{\tau}^{(a,b)}<\infty]
\nonumber\\
&=\mathbb{E}_y[e^{-\delta\widetilde{\tau}^{-}_{0} }~;~\widetilde{\tau}^{-}_{0} <\infty]\nonumber\\
&=\mathbb{J}_{a}^{(\delta,\gamma)}(y-a)
-\mathbb{I}_{a}^{(\delta,\gamma)}(y-a)
\frac{\delta\mathbb{W}^{(\delta)}(a)\mathbb{Z}^{(\delta)}(a,\phi(\delta+\gamma))}{\mathbb{Z}^{(\delta)'}(a,\phi(\delta+\gamma))},
\end{align}
where the last equation is obtained from equation $(\ref{3.4})$ in Lemma $\ref{lem3.2}$.

For the case $a\le y<b$, according to the law of total probability and the strong Markov property, we have%
\begin{align}\label{15}
&L(y;a,b)\\&=\mathbb{E}_y[e^{-\delta\widetilde{\tau}^{(a,b)}}~;~\widetilde{\tau}^{(a,b)}<\infty]\nonumber \\
&=\mathbb{E}_y[e^{-\delta\widetilde{\tau}^{(a,b)}}~;~\widetilde{\tau}^{(a,b)}<\infty,~
\tau^{-}_{0} <\tau^{+}_{b} ]+\mathbb{E}_y[e^{-\delta\widetilde{\tau}^{(a,b)}}~;~\widetilde{\tau}^{(a,b)}<\infty,~
\tau^{-}_{0} >\tau^{+}_{b} ]\nonumber\\
&=\mathbb{E}_y[e^{-\delta\tau^{-}_{0} }~;~\tau^{-}_{0} <\infty,~
\tau^{-}_{0} <\tau^{+}_{b} ]+\mathbb{E}_y[e^{-\delta\tau^{+}_{b} }\mathbb{E}_b[e^{-\delta\widetilde{\tau}^{(a,b)}}~;~\widetilde{\tau}^{(a,b)}<\infty]
~;
\tau^{+}_{b} <\tau^{-}_{0} ]\nonumber\\
&=\mathbb{E}_y[e^{-\delta\tau^{-}_{0} }~;
\tau^{-}_{0} <\tau^{+}_{b} ]
+\mathbb{E}_y[e^{-\delta\tau^{+}_{b} }~;
\tau^{+}_{b} <\tau^{-}_{0} ]L(b;a,b)\nonumber\\
&
=\left(Z^{(\delta)}(y)-Z^{(\delta)}(b)\frac{W^{(\delta)}(y)
}{W^{(\delta)}(b)}\right)+
\frac{W^{(\delta)}(y)}{W^{(\delta)}(b)}
L(b;a,b),
\end{align}
where $L(b;a,b)=\mathbb{J}_{a}^{(\delta,\gamma)}(b-a)
-\mathbb{I}_{a}^{(\delta,\gamma)}(b-a)
\frac{\delta\mathbb{W}^{(\delta)}(a)\mathbb{Z}^{(\delta)}
(a,\phi(\delta+\gamma))}{\mathbb{Z}^{(\delta)'}(a,\phi(\delta+\gamma))}
$
is obtained by equation $(\ref{14})$ when $y=b$ and the last equation is obtained by using Lemma $\ref{lem3.3}$.

Similarly, for the last case $0\le y<a$, we have%
\begin{align}\label{16}
&L(y;a,b)\\&=\mathbb{E}_y[e^{-\delta\widetilde{\tau}^{(a,b)}}~;~\widetilde{\tau}^{(a,b)}<\infty]\nonumber \\
&=\mathbb{E}_y[e^{-\delta\widetilde{\tau}^{(a,b)}}~;~\widetilde{\tau}^{(a,b)}<\infty,~
\tau^{-}_{0} <\tau^{+}_{a} ]+\mathbb{E}_y[e^{-\delta\widetilde{\tau}^{(a,b)}}~;~\widetilde{\tau}^{(a,b)}<\infty,~
\tau^{-}_{0} >\tau^{+}_{a} ]\nonumber\\
&=\mathbb{E}_y[e^{-\delta\tau^{-}_{0} }~;~\tau^{-}_{0} <\infty,~
\tau^{-}_{0} <\tau^{+}_{a} ]+\mathbb{E}_y[e^{-\delta\tau^{+}_{a} }\mathbb{E}_a[e^{-\delta\widetilde{\tau}^{(a,b)}}~;~\widetilde{\tau}^{(a,b)}<\infty]
~;
\tau^{+}_{a} <\tau^{-}_{0} ]\nonumber\\
&=\mathbb{E}_y[e^{-\delta\tau^{-}_{0} }~;
\tau^{-}_{0} <\tau^{+}_{a} ]
+\mathbb{E}_y[e^{-\delta\tau^{+}_{a} }~;
\tau^{+}_{a} <\tau^{-}_{0} ]L(a;a,b)\nonumber\\
&
=\left(Z^{(\delta)}(y)-Z^{(\delta)}(a)\frac{W^{(\delta)}(y)}{W^{(\delta)}(a)}\right)+
\frac{W^{(\delta)}(y)}{W^{(\delta)}(a)}L(a;a,b),
\end{align}
where $L(a;a,b)=Z^{(\delta)}(a)+\frac{W^{(\delta)}(a)}{W^{(\delta)}(b)}
\left\{\mathbb{J}_{a}^{(\delta,\gamma)}(b-a)
-Z^{(\delta)}(b)
-\mathbb{I}_{a}^{(\delta,\gamma)}(b-a)
\frac{\delta\mathbb{W}^{(\delta)}(a)\mathbb{Z}^{(\delta)}(a,\phi(\delta+\gamma))}{\mathbb{Z}^{(\delta)'}(a,\phi(\delta+\gamma))}
\right\}$ is obtained by $y=b$ in equation $(\ref{15})$. This proof is end. $\Box$

\begin{rmk}\label{rmk.2}~~{\rm
 Noting that by taking $\delta\rightarrow 0$ in equations $(\ref{14})$-$(\ref{16})$, we obtain that
 for $b\ge a>0,~y\ge 0$, $\widetilde{\tau}^{(a,b)}<\infty$ $\mathbb{P}_y$$-a.s.$ as $\lim_{\delta\rightarrow 0}L(y;a,b)=1$.
}
\end{rmk}

Let $L^R(y;b)$, $y\ge 0$, denote the Laplace transform of the ruin time for the spectrally negative L\'{e}vy  processes with the ratcheting dividend strategy(single-rise). By taking $\gamma\downarrow0$ in Theorems $\ref{thm4.1}$, we derive the expression of Laplace transform $L^R(y;b)$, which can reduce the results in Section $5$ of Albrecher, B\"{a}uerle, and Bladt (2018).
\begin{cory}\label{cory.3}~~For $\delta, c_1, y\ge 0$ and $c_2, b>0$, we have
\begin{align*}
L^R(y;b)=
\left\{
  \begin{aligned}
  \mathbb{Z}^{(\delta)}(y)-\frac{\delta}{ \phi(\delta) }\mathbb{W}^{(\delta)}(y)&,~~y\ge b,\\
   Z^{(\delta)}(y)+
  \frac{W^{(\delta)}(y)}{W^{(\delta)}(b)}\left( \mathbb{Z}^{(\delta)}(b)-Z^{(\delta)}(b)-\frac{\delta}{ \phi(\delta) }\mathbb{W}^{(\delta)}(b) \right)&,~~0\le y< b.
  \end{aligned}
  \right.
\end{align*}
\end{cory}


\section{Numerical illustrations}\label{Num}
In this subsection, we give some graphs, under the special cases of the Brownian motion with drift, to show the effects of some parameters on the expected NPV of dividends up to ruin and the Laplace transform of the ruin time.

Let $U(t)=y+\mu t+\sigma B(t),~t\ge 0,$
where $y>0$ is the initial surplus, $\mu>0$ is a constant drift, $\sigma>0$ and $\{B(t),~t\ge0\}$ denotes a
standard Brownian motion. The Brownian risk model has also been considered by Albrecher, B\"{a}uerle, and Bladt (2018). It is worth mentioning that the expected gain per unit time should be positive( i.e. $\mu>c_1+c_2$).
In the diffusion case the scale functions of $X$ and $\widetilde{X}$ are given by
\begin{align*}
&W^{(\delta)}(x)=\kappa( e^{\theta_1 x}-e^{\theta_2 x}),~~~~
W^{(\delta+\gamma)}(x)=l( e^{\theta_3 x}-e^{\theta_4 x}),\\
&\mathbb{W}^{(\delta)}(x)=\widetilde{\kappa}( e^{\widetilde{\theta_1} x}-e^{\widetilde{\theta_2} x}),~~~~
\mathbb{W}^{(\delta+\gamma)}(x)=\widetilde{l}( e^{\widetilde{\theta_3} x}-e^{\widetilde{\theta_4} x}),
\end{align*}
where
\begin{align*}
&\kappa:= ((\mu - c_1)^2 + 2 \sigma^2 \delta)^{ -1/2 },\\
&\widetilde{\kappa}:= ((\mu - c_1 - c_2)^2 + 2 \sigma^2 \delta)^{-1/2},\\
&l:= ((\mu - c_1)^2 + 2 \sigma^2 (\delta+\gamma))^{ -1/2 },\\
&\widetilde{l}:= ((\mu - c_1 - c_2)^2 + 2 \sigma^2 (\delta+\gamma))^{-1/2},\\
&\theta_1=\frac{-(\mu-c_1)+ \sqrt{(\mu-c_1)^2+2\delta\sigma^2}}{\sigma^2},\\
&\theta_2=\frac{-(\mu-c_1 )- \sqrt{(\mu-c_1 )^2+2\delta\sigma^2}}{\sigma^2},\\
&\widetilde{\theta}_1=\frac{-(\mu-c_1-c_2)+ \sqrt{(\mu-c_1-c_2)^2+2\delta\sigma^2}}{\sigma^2},\\
&\widetilde{\theta}_2=\frac{-(\mu-c_1-c_2)- \sqrt{(\mu-c_1-c_2)^2+2\delta\sigma^2}}{\sigma^2},\\
&\theta_3=\frac{-(\mu-c_1)+ \sqrt{(\mu-c_1)^2+2(\delta+\gamma)\sigma^2}}{\sigma^2},\\
&\theta_4=\frac{-(\mu-c_1 )- \sqrt{(\mu-c_1 )^2+2(\delta+\gamma)\sigma^2}}{\sigma^2},\\
&\widetilde{\theta}_3=\frac{-(\mu-c_1-c_2)+ \sqrt{(\mu-c_1-c_2)^2+2(\delta+\gamma)\sigma^2}}{\sigma^2},\\
&\widetilde{\theta}_4=\frac{-(\mu-c_1-c_2)- \sqrt{(\mu-c_1-c_2)^2+2(\delta+\gamma)\sigma^2}}{\sigma^2}.
\end{align*}
Note that, in this case $\phi(\delta+\gamma)=\widetilde{\theta}_3$.
For $x,\delta\ge 0$, by the define of $Z^{(\delta)}(x)$ and $\mathbb{Z}^{(\delta)}(x)$, we have
\begin{align*}
Z^{(\delta)}(x)=\frac{\delta}{\theta_1} {W}^{(\delta)}(x) + e^{\theta_2 x},~~~~
\mathbb{Z}^{(\delta)}(x)=\frac{\delta}{\widetilde{\theta}_1} \mathbb{W}^{(\delta)}(x) + e^{\widetilde{\theta}_2 x}.
\end{align*}
By some algebraic manipulations, for $y\ge a\ge 0$, we have
\begin{align*}
&\mathbb{Z}^{(\delta)}(a,\phi(\delta+\gamma))=-\frac{\gamma}{\widetilde{\theta}_1-\widetilde{\theta}_3}
 \mathbb{W}^{(\delta)}(a)+e^{\widetilde{\theta}_2 a },~~\\
&\mathbb{Z}^{(\delta)'}(a,\phi(\delta+\gamma))=\frac{\gamma \widetilde{\theta}_1}{
\widetilde{\theta}_1 - \widetilde{\theta}_3}
\mathbb{W}^{(\delta)}(a) + \widetilde{\theta}_3 e^{\widetilde{\theta}_2 a },\\
&\overline{\overline{{W}}}^{(\delta+\gamma)}(y-a)=\frac{1}{{\theta}_3^2}{W}^{(\delta+\gamma)}(y-a)
+\frac{\mu - c_1}{(\delta+\gamma)^2} e^{{\theta}_4 (y-a)}
- \frac{\mu - c_1}{(\delta+\gamma)^2}
- \frac{y-a}{\delta+\gamma},\\
&\overline{\overline{\mathbb{W}}}^{(\delta+\gamma)}(y-a)=\frac{1}{\widetilde{\theta}_3^2}\mathbb{W}^{(\delta+\gamma)}(y-a)
+\frac{\mu - c_1-c_2}{(\delta+\gamma)^2} e^{\widetilde{\theta}_4 (y-a)}
- \frac{\mu - c_1-c_2}{(\delta+\gamma)^2}
- \frac{y-a}{\delta+\gamma},\\
&{I}_{a}^{(\delta,\gamma)}(y-a)=
\frac{\gamma(2\theta_3-\theta_1)}{\theta_3(\theta_1-\theta_3)}W^{(\delta+\gamma)}(y-a)
 + \frac{\delta }{\delta+\gamma} e^{\theta_4 (y-a)}
+ \frac{e^{\theta_2 (a)}}{{W}^{(\delta)}(a) } W^{(\delta+\gamma)}(y-a)
+\frac{1}{\delta+\gamma},\\
&\mathbb{I}_{a}^{(\delta,\gamma)}(y-a)=
\frac{\gamma(2\widetilde{\theta}_3-\widetilde{\theta}_1)}{\widetilde{\theta}_3(\widetilde{\theta}_1-\widetilde{\theta}_3)}\mathbb{W}^{(\delta+\gamma)}(y-a)
 + \frac{\delta }{\delta+\gamma} e^{\widetilde{\theta}_4 (y-a)}
+ \frac{e^{\widetilde{\theta}_2 (a)}}{\mathbb{W}^{(\delta)}(a) } \mathbb{W}^{(\delta+\gamma)}(y-a)
+\frac{1}{\delta+\gamma},\\
&{{J}_{a}^{(\delta,\gamma)}(y-a)}=
 (1 +\frac{l}{ k } )e^{\theta_2 y}
-\frac{\gamma\delta }{ \theta_3(\theta_1 - \theta_3)}  W^{(\delta)}(a) W^{(\delta+\gamma)}(y-a)
+\frac{\delta^2}{\theta_1(\delta+\gamma)} W^{(\delta)}(a) e^{\theta_4(y-a)}\\
&~~~~~~~~~~~~~~~~ +\frac{\gamma\delta}{\theta_1(\delta+\gamma)} W^{(\delta)}(a)
-( \frac{\gamma}{\theta_3}-\frac{\delta }{\theta_1}+ \frac{\gamma}{\theta_2 - \theta_3} )
  W^{(\delta+\gamma)}(y-a)e^{ \theta_2 a}
\\&~~~~~~~~~~~~~~~~- ( \frac{\gamma}{\delta+\gamma}+ \frac{l}{ k }  )
  e^{\theta_4(y-a)+\theta_2 a}
+ \frac{\gamma}{\delta+\gamma}e^{\theta_2 a},\\
&{\mathbb{J}_{a}^{(\delta,\gamma)}(y-a)}=
 (1 +\frac{\widetilde{l}}{ \widetilde{k} } )e^{\widetilde{\theta}_2 y}
-\frac{\gamma\delta }{ \widetilde{\theta}_3(\widetilde{\theta}_1 - \widetilde{\theta}_3)}  \mathbb{W}^{(\delta)}(a) \mathbb{W}^{(\delta+\gamma)}(y-a)
+\frac{\delta^2}{\widetilde{\theta}_1(\delta+\gamma)} \mathbb{W}^{(\delta)}(a) e^{\widetilde{\theta}_4(y-a)}\\
&~~~~~~~~~~~~~~~~ +\frac{\gamma\delta}{\widetilde{\theta}_1(\delta+\gamma)} \mathbb{W}^{(\delta)}(a)
-( \frac{\gamma}{\widetilde{\theta}_3}-\frac{\delta }{\widetilde{\theta}_1}+ \frac{\gamma}{\widetilde{\theta}_2 - \widetilde{\theta}_3} )
  \mathbb{W}^{(\delta+\gamma)}(y-a)e^{ \widetilde{\theta}_2 a}
\\
&~~~~~~~~~~~~~~~~- ( \frac{\gamma}{\delta+\gamma}+ \frac{l}{ k }  )
  e^{\widetilde{\theta}_4(y-a)+\widetilde{\theta}_2 a}
+ \frac{\gamma}{\delta+\gamma}e^{\widetilde{\theta}_2 a}.
\end{align*}
\subsection{Analysis with the expected NPV of dividends up to ruin}\label{ANA_div}
In this section, we reveal the impact of various parameters on the expected discounted cumulative dividend function. In order to investigate that, in the following analysis, unless otherwise specified, the basic parameter settings are as follows: $\mu=1$, $\sigma=2$, $c_1=0$, $c_2=0.1$, $\gamma=1$, $\delta=0.05$.

\begin{figure}[htbp]
\centering
\subfigure[]{
\includegraphics[scale=0.35]{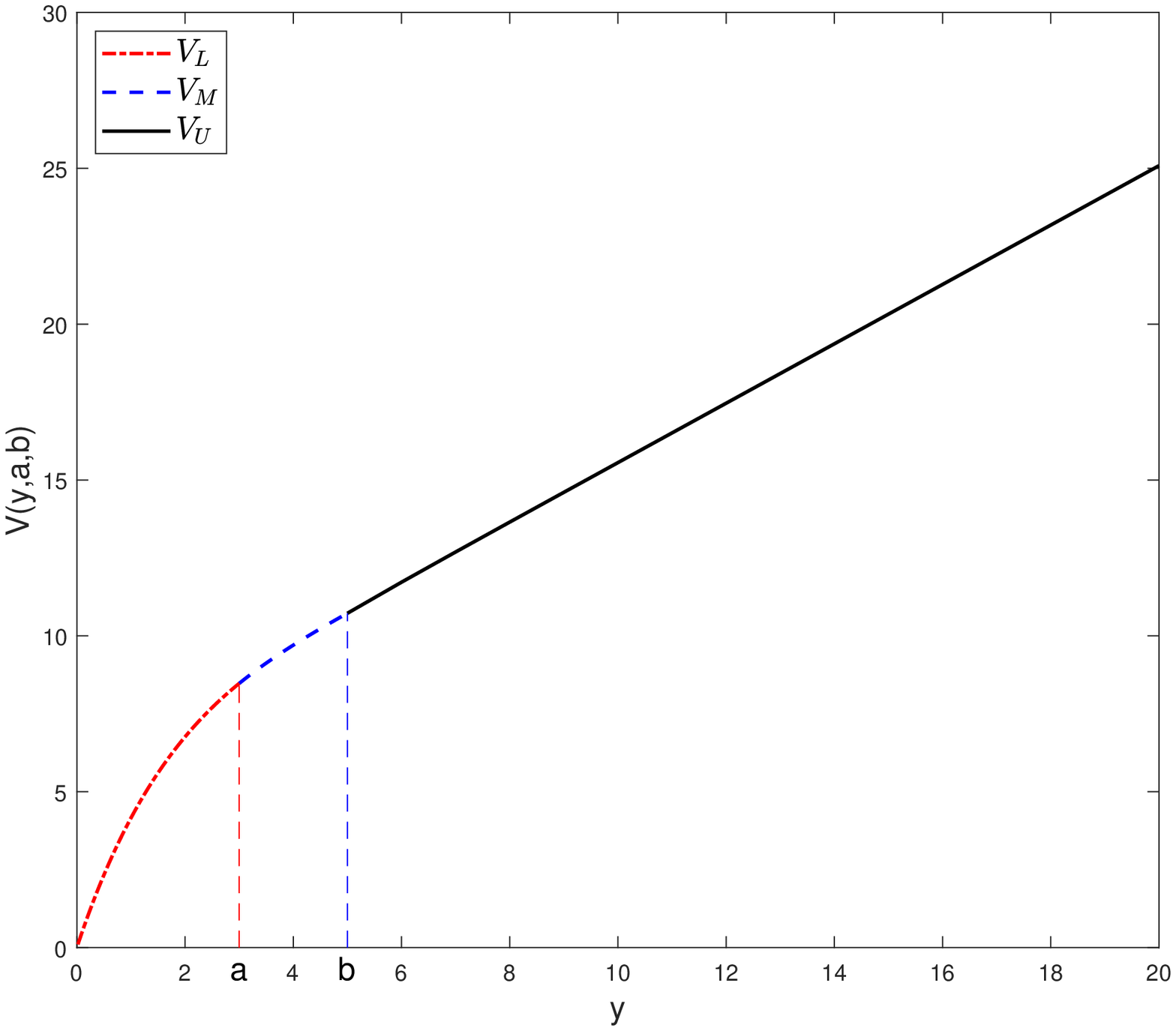}
}%
\subfigure[]{
\includegraphics[scale=0.35]{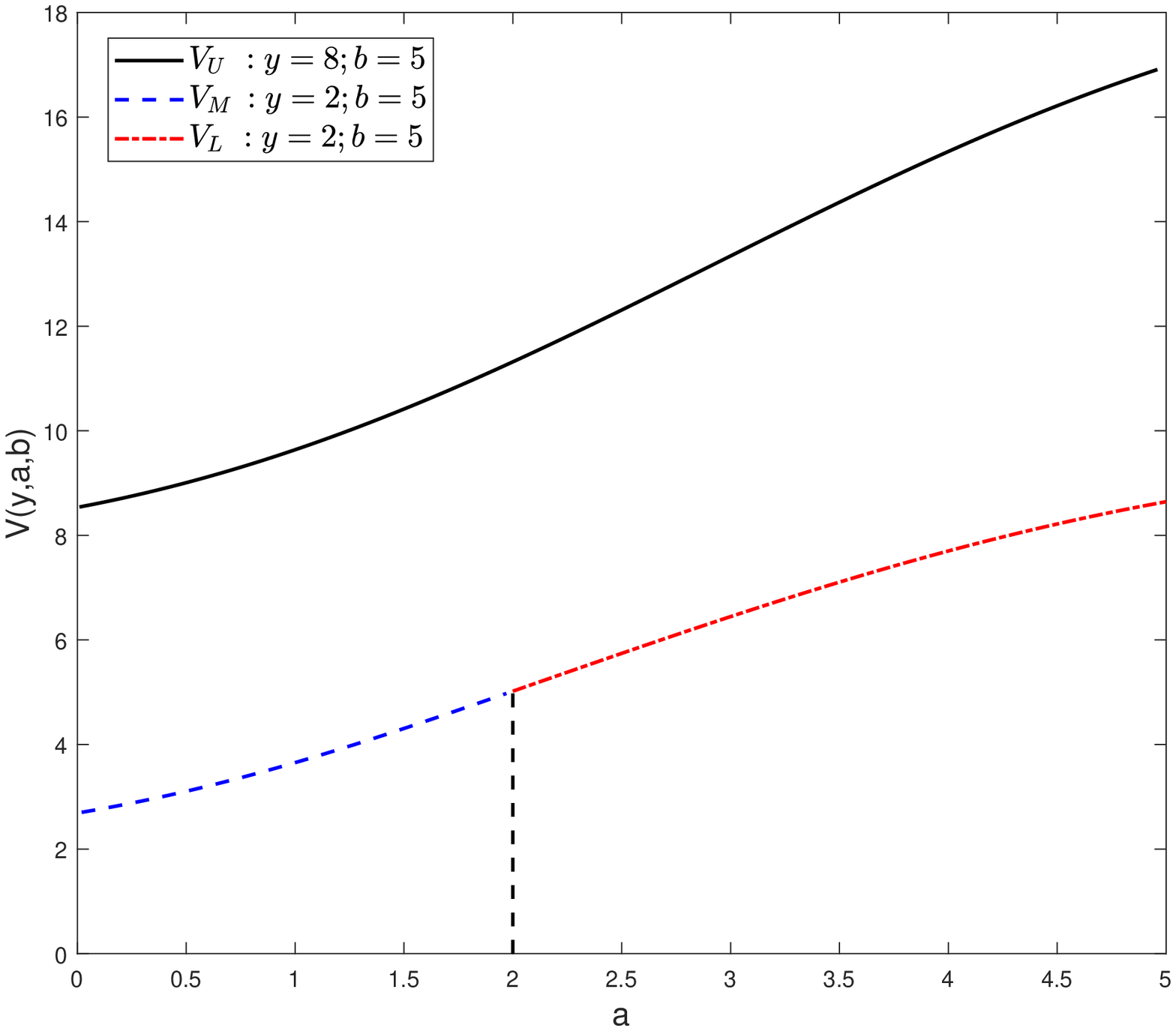}
}%
\caption{ The expected net present value of dividends paid up to ruin $V(y;a,b)$ (\textbf{a}): as a function of
$y$ for three different intervals $y\in[0,a)$(dotted and dashed curve), $y\in [a,b)$(dashed
curve), and $y\in [b,\infty)$(solid curve); 
(\textbf{b}): as a function of
$a$ when $y=8;b=5$(solid), $y=2;b=5$(dashed), and $y=2;b=5$(dotted and dashed), respectively.
}
\label{fig.1}
\end{figure}

\begin{figure}[htbp]
\centering
\subfigure[]{
\includegraphics[scale=0.35]{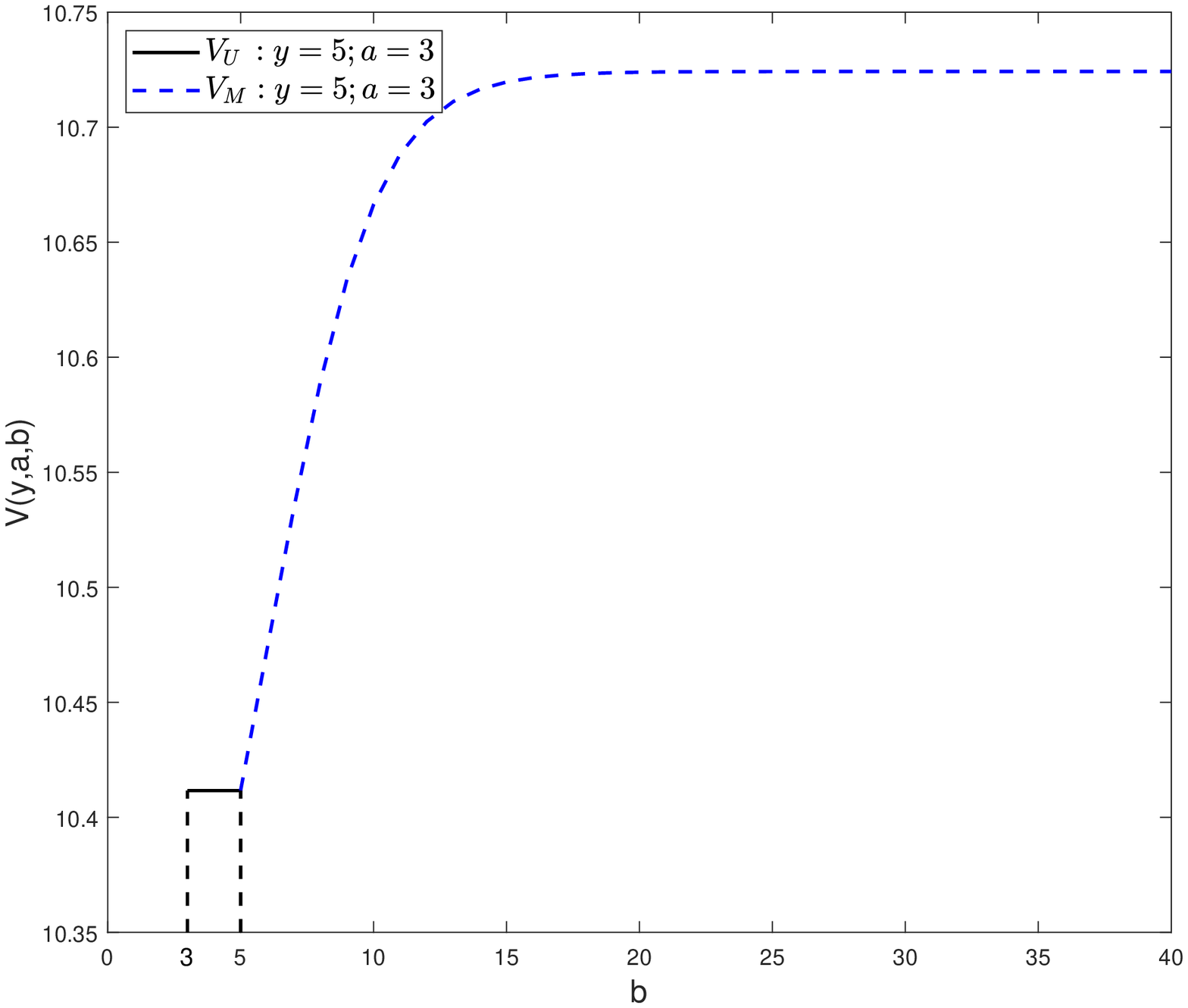}
}%
\subfigure[]{
\includegraphics[scale=0.35]{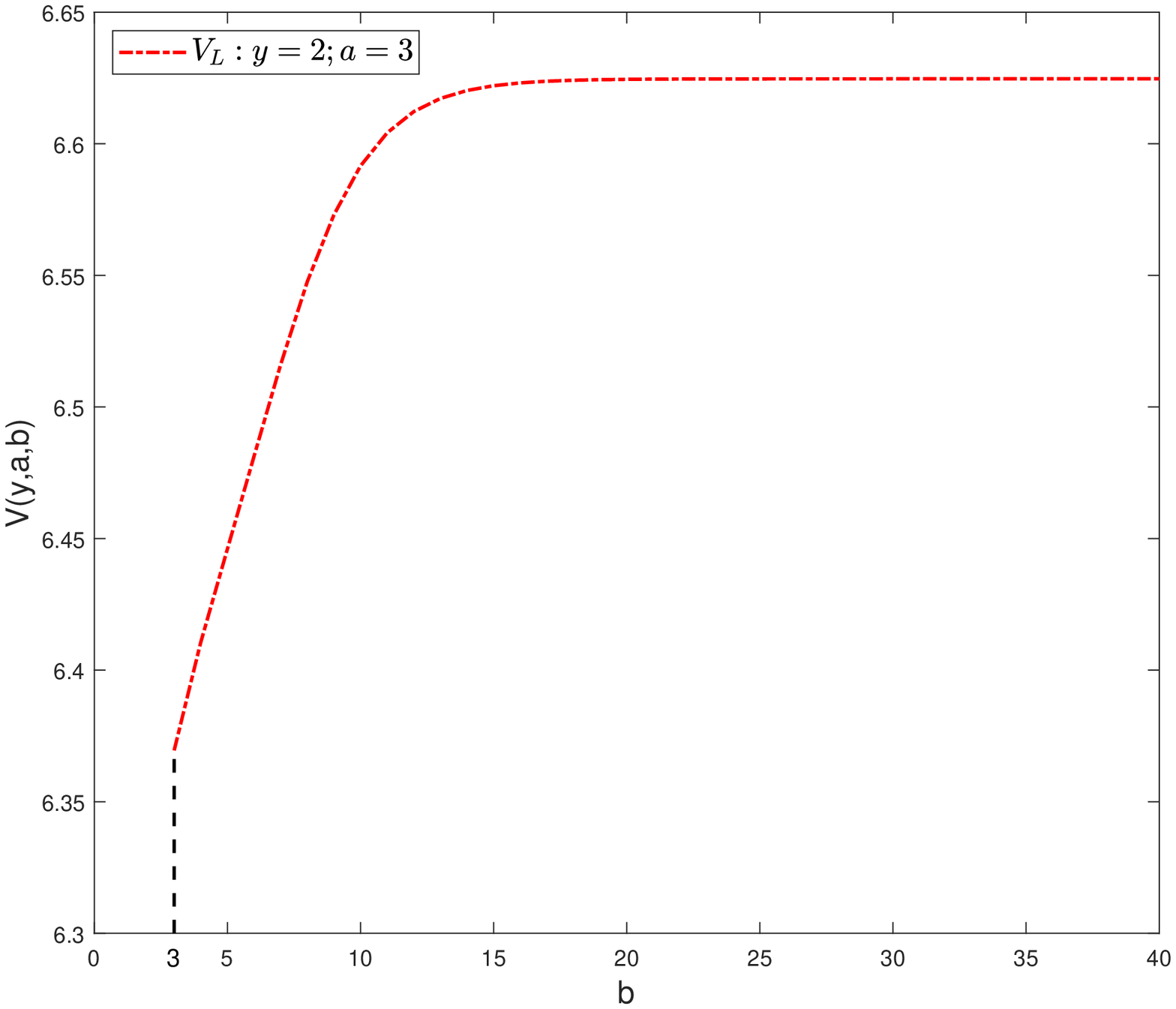}
}
\caption{The expected net present value of dividends paid up to ruin $V(y;a,b)$ as a function of
$b$  when $y=5;a=3$(solid in left), $y=5;a=3$(dashed in left), and $y=2;a=3$(dotted and dashed in right), respectively.
 }
\label{fig.3}
\end{figure}

\begin{figure}[hbp]
\centering
\subfigure[]{
\includegraphics[scale=0.35]{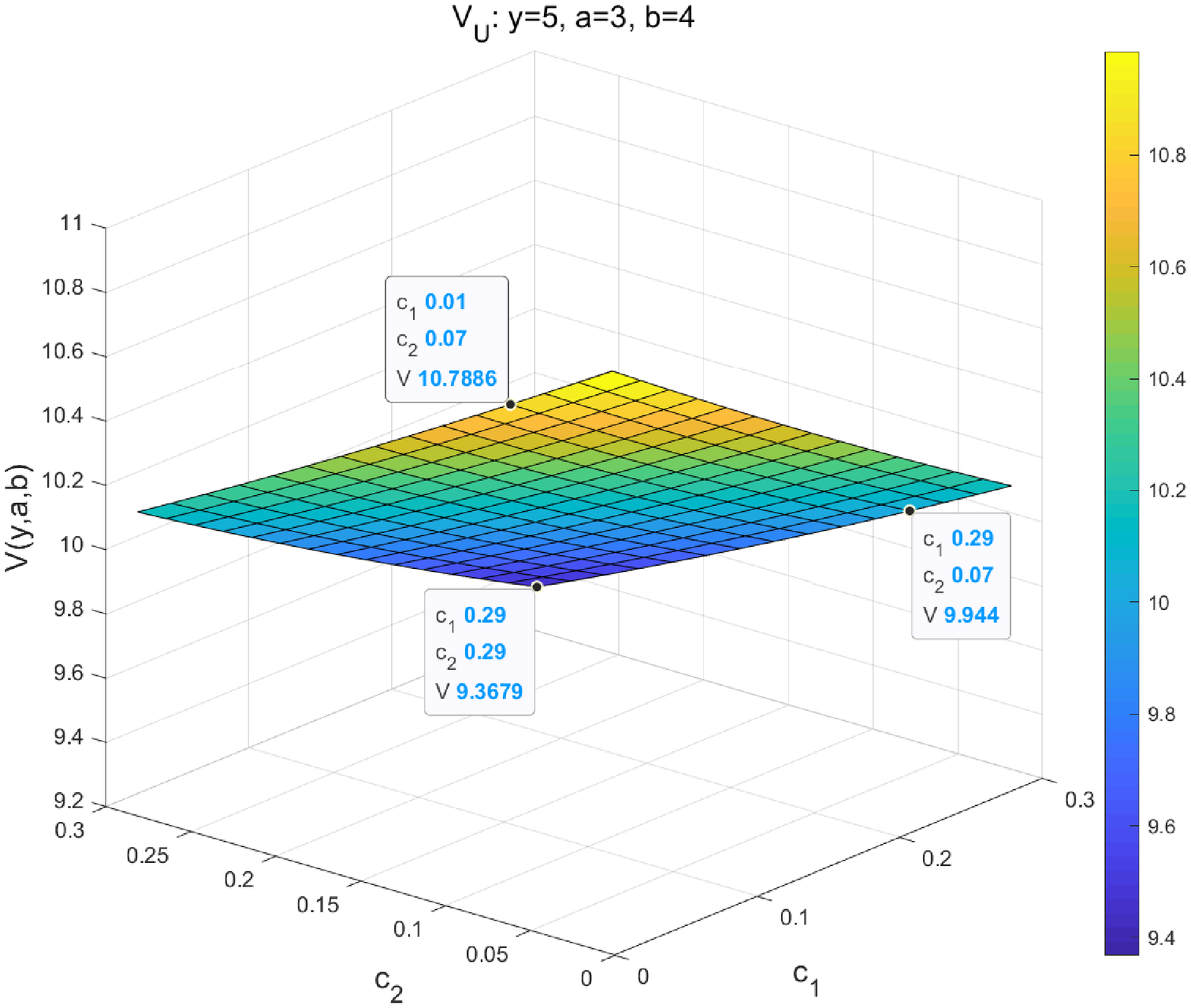}
}%
\subfigure[]{
\includegraphics[scale=0.35]{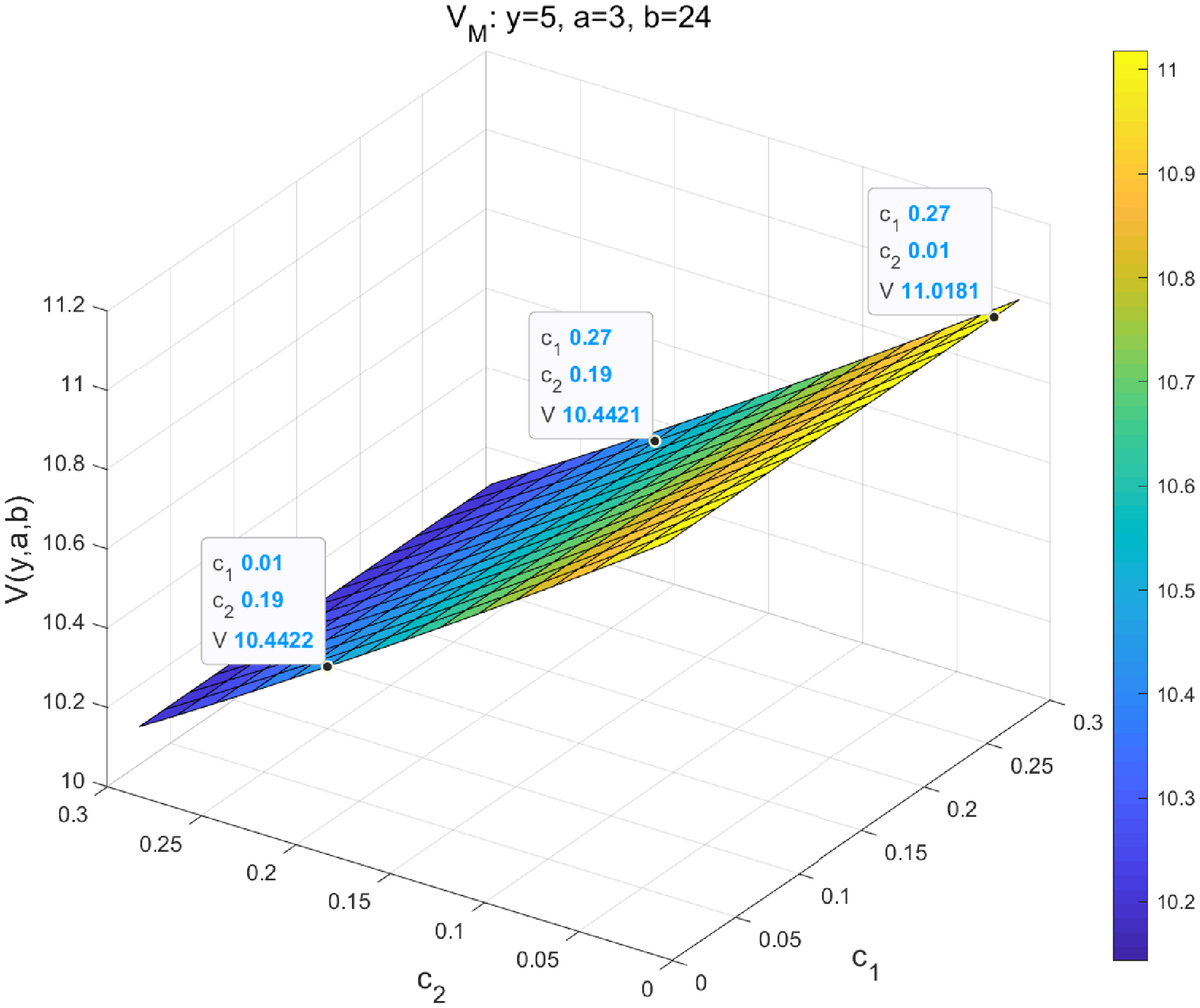}
}%

\subfigure[]{
\includegraphics[scale=0.35]{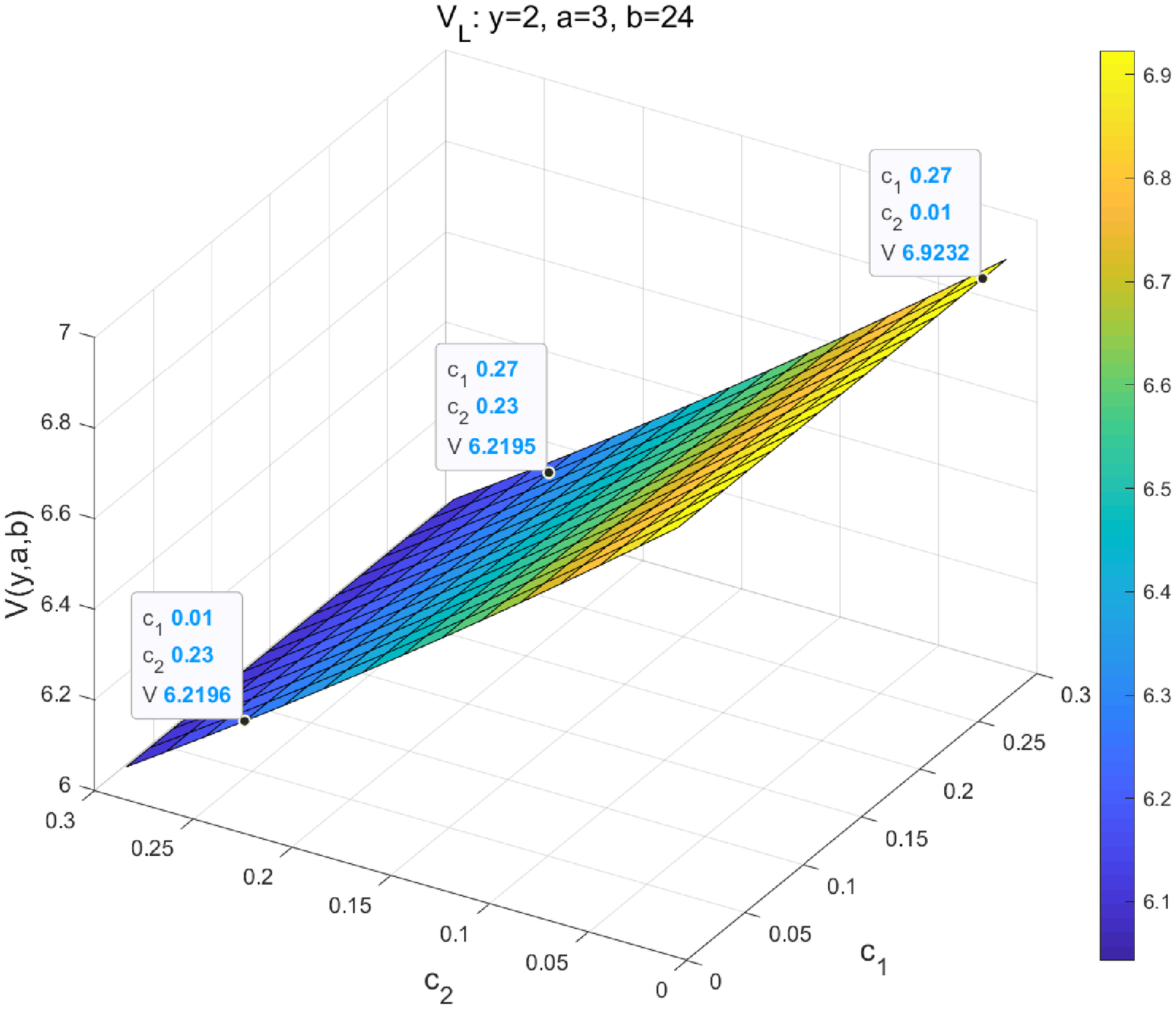}
}%
\caption{The expected net present value of dividends paid up to ruin $V(y;a,b)$ as a function of
$(c_1, c_2)$ for $(c_1, c_2)$$\in[0,0.3]\times[0,0.3]$,
for three different cases£º
(\textbf{a}) $y=5;a=3;b=4; $~
(\textbf{b}) $y=5;a=3;b=24; $
(\textbf{c}) $y=2;a=3;b=24. $~
}
\label{fig.4}
\end{figure}

\begin{figure}[hp]
\centering
\includegraphics[scale=0.65]{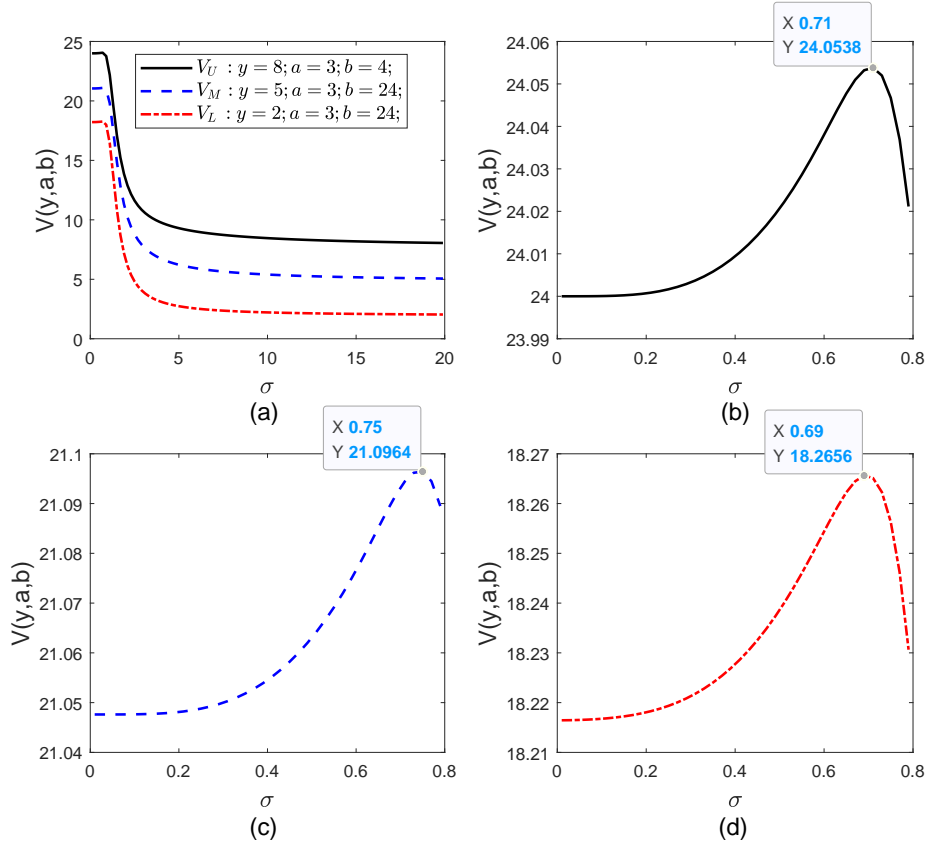}
\caption{The expected net present value of dividends paid up to ruin $V(y;a,b)$ as a function of
$\sigma$ for $(c_1, c_2)$$\in[0,0.25]\times[0,0.25]$,
for three different cases $y=8;a=3;b=4 $(solid), $y=5;a=3;b=24 $(dashed),~and~$y=2;a=3;b=24 $(dotted and dashed). In order to easily see the trend of $V(y;a,b)$ with respect to $\sigma$, we add sub-figures (\textbf{b}-\textbf{d}).
}
\label{fig.5}
\end{figure}

For finding the influence, in Figures $\ref{fig.1}$-$\ref{fig.5}$, we show the curves of $V(y;a,b)$ as functions of $y $, $a $, $b $, $ (c_1, c_2)$, and $ \sigma$, respectively. First, for fixed ratcheting-periodic barriers $a=3$ and $b=5$, we find from Figure $\ref{fig.1}$(a) that the expected NPV of dividends paid up to ruin $V(y;a,b)$ increases as the initial value $y$ increases, which is an intuitive conclusion. When the initial surplus increases, the surplus process is more likely to be above the ratcheting-periodic barriers, and is more likely to survive longer, and then there may be more dividend payments paid off. Meanwhile, we see that $V(y;a,b)$ is continuous obviously in $y$, which coincides with Remark $\ref{rmk.1}$.
Figure $\ref{fig.1}$(b) shows a result that be not intuitively understood: For fixed $y=2 ~or~8$ and $b=5$, $V(y;a,b)$ increases with respect to the periodic dividend barrier $a$, which is also a very important discovery.
We can explain this phenomenon as follows: As $a$ increases, the ruin time may be delayed, which result that the potential dividend at later times increases.
From that we also find the optimal periodic value $a^*=b$. This implies that under certain circumstances, in order to maximize $V(y;a,b)$, we should let the periodic dividend barrier be equal to the ratcheting dividend barrier(i.e. $a=b$). In fact, Song and Sun (2021) have been studied this situation in a dual risk model. In view of $a^*=b$, we want to find the optimal ratcheting value $b^*$, which is corresponding to the maximum value of $V(y;a,b)$ for fixed the values of other parameters except parameter $b$. Then we plot curves for $V(y;a,b)$ as a function of the ratcheting barrier $b$ in Figure $\ref{fig.3}$.

From Figure $\ref{fig.3}(a)$, we see that the value of $V(y;a,b)$ is a constant when the initial value $y$ is greater than $b$. This is consistent with Theorem $\ref{thm3.1}$: the value of $V_U$ has nothing to do with $b$. Figure $\ref{fig.3}$ also shows that: $V(y;a,b)$ is increasing in the ratcheting dividend barrier $b$ when $y$ is lower than $b$ and converges to a nonzero constant when $b$ tends to infinity. Then the optimal periodic dividend barrier $b^*=\infty$. However it is unrealistic to set $b^*=\infty$ for nondecreasing dividend rate, so in the light of the result of Figure $\ref{fig.3}$, we can find the approximate value $\widetilde{b^*}$ of optimal ratcheting barrier within an acceptable error range.

As can be seen in Figure $\ref{fig.4}$, for three different cases, all of $V(y;a,b)$ are decreasing functions of $c_1$ and $c_2$, respectively. Meanwhile, we find from the remarked values of $V(y;a,b)$ that the influence of $c_1$ and $c_2$ on $V(y;a,b)$ is relatively smaller. For this result, we think the optimal barriers $c_1^*$ and $c_2^*$ do not make sense. Conversely, we can choose the value of $c_1$ and $c_2$ depending on the market condition and the insurance company's own needs.

Next, we describe the overall trend of $V(y;a,b)$ as the market volatility $\sigma$ changes in Figure $\ref{fig.5}$(a). On one hand, we find that $V(y;a,b)$ is not a monotone function of $\sigma$ and converges to a nonzero constant as $\sigma\rightarrow \infty$. This shows that even when the market volatility is very high, $V(y;a,b)$ can remain at a relatively stable value. On the other hand, for making a more precise conclusion, we plot $V(y;a,b)$ as functions of $\sigma\in(0,0.8)$ in Figure $\ref{fig.5}$(b-d). From that we see for every fixed $y$, $a$, and $b$, that $V(y;a,b)$ is a concave function of $\sigma$, which means that $V(y;a,b)$ first increases and then decreases in $\sigma$. From the remarked values in Figure $\ref{fig.5}$(b-d)(i.e. the corresponding optimal value $V(y;a,b)|_{\sigma=\sigma^*}$), the value $\sigma^*$ under these three cases is different. This indicates that $\sigma^*$ depends on other parameters.
Therefore, we can control other controllable parameters to make $\sigma$ close to $\sigma^*$, so as to obtain the maximum value of $V(y;a,b)$.
Finally, as $\sigma\rightarrow 0$ the process becomes deterministic, but $\sigma^*\neq 0$, which interprets that no turbulent market is not optimal and a smaller market volatility can increase $V(y;a,b)$.

\subsection{Analysis with the Laplace transform of the ruin time}\label{ANA_ruin}

\begin{figure}[hp]
\centering
\subfigure[]{
\includegraphics[scale=0.3]{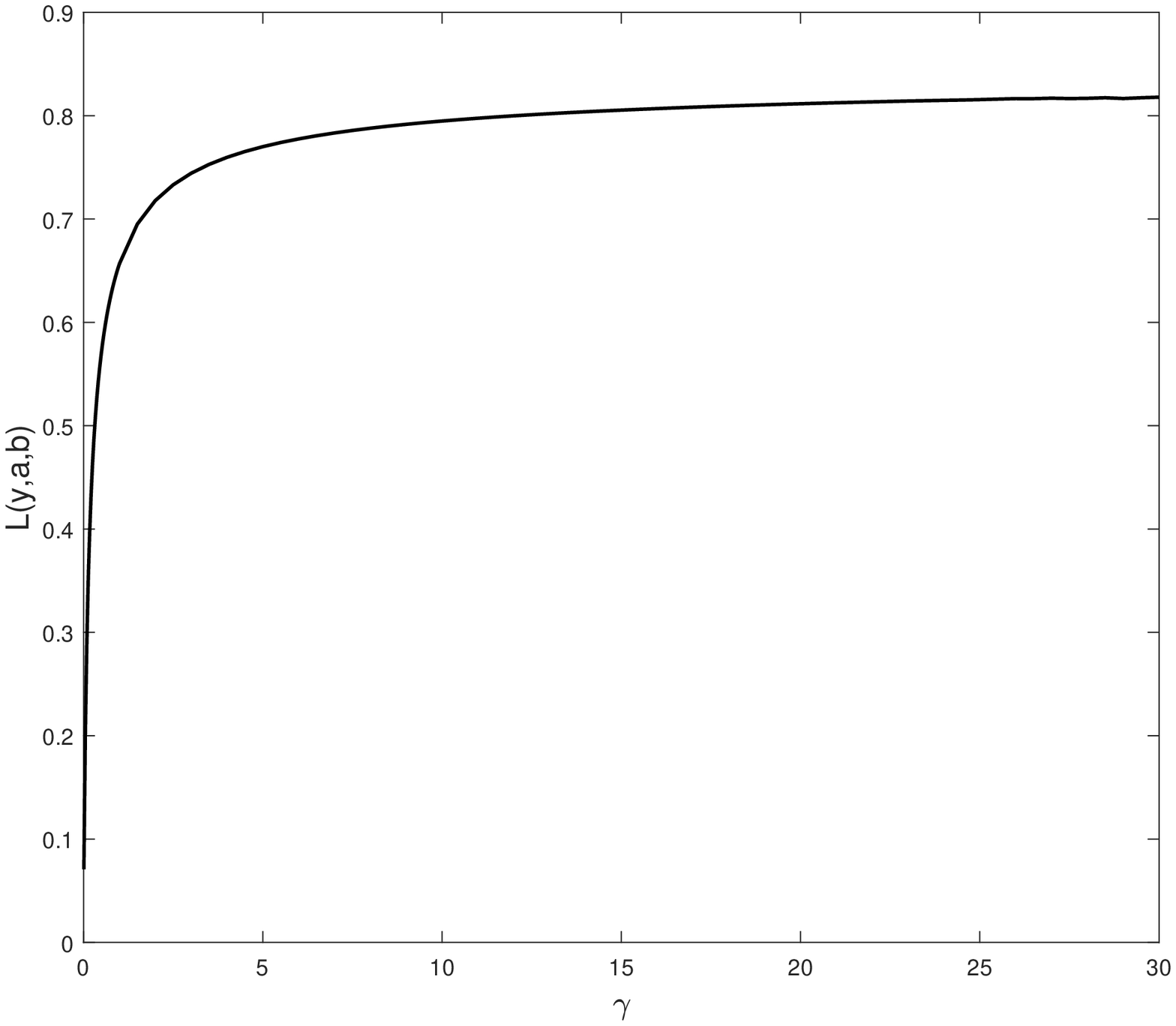}
}%
\subfigure[]{
\includegraphics[scale=0.3]{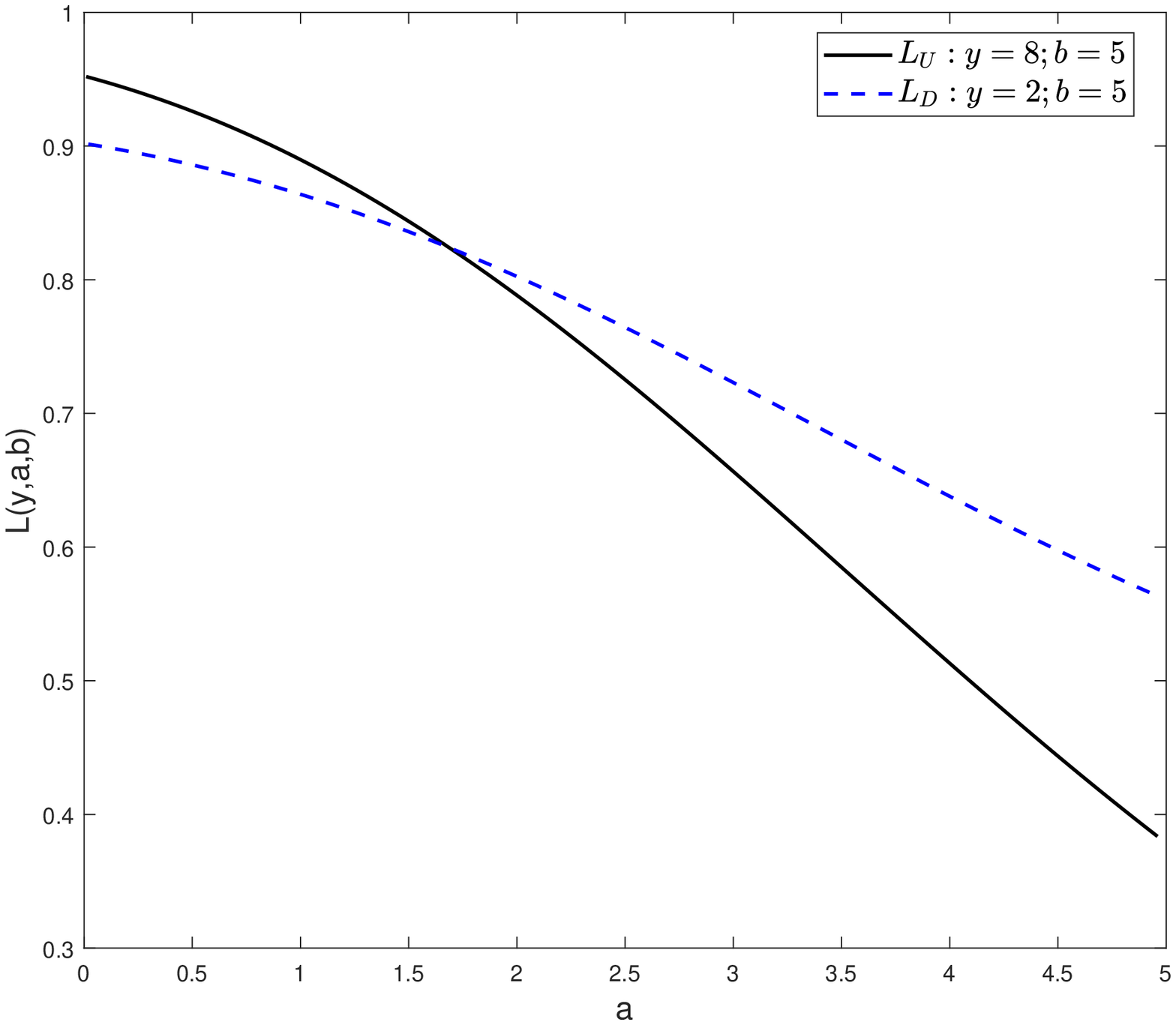}
}%

\subfigure[]{
\includegraphics[scale=0.3]{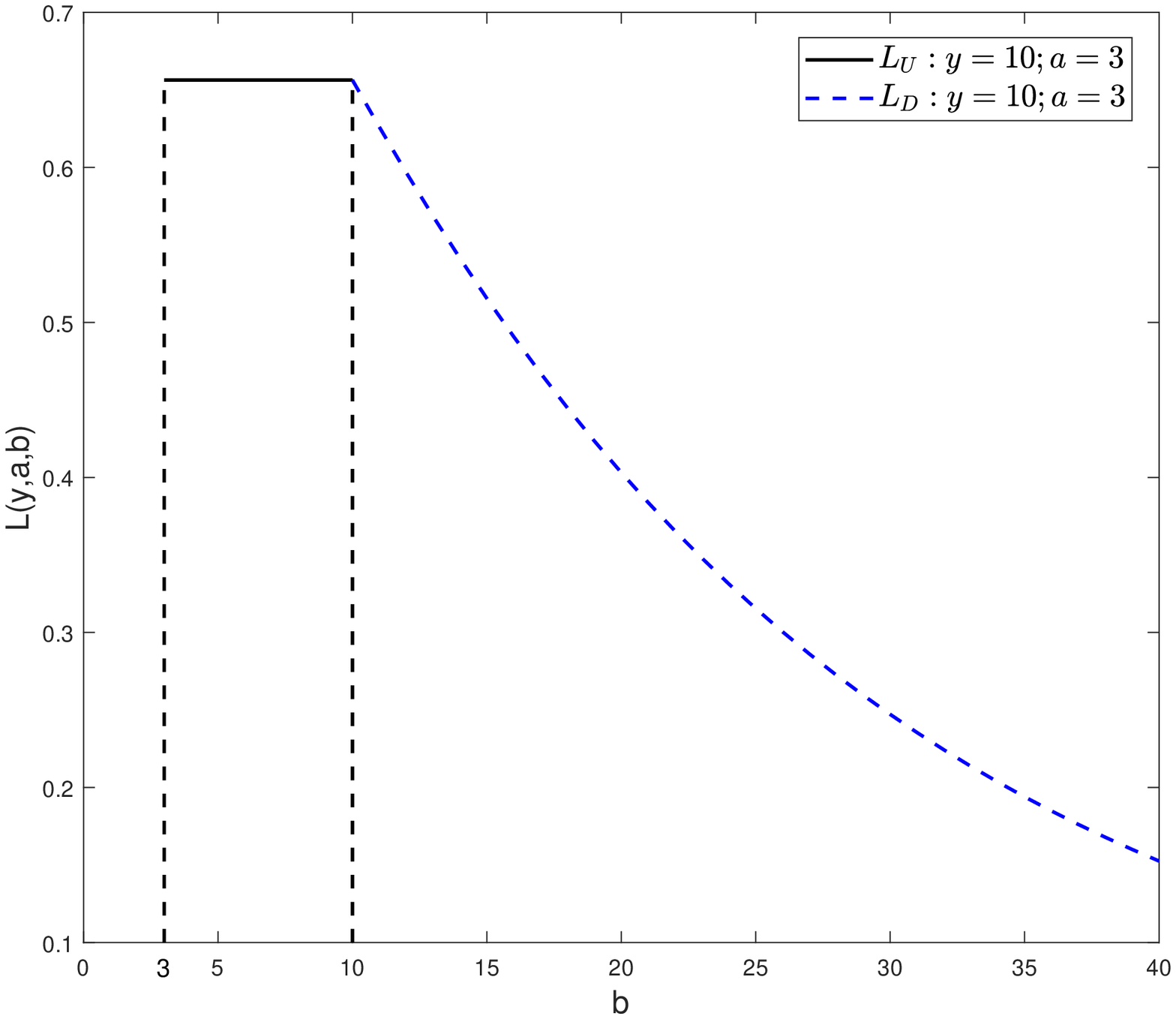}
}%
\subfigure[]{
\includegraphics[scale=0.3]{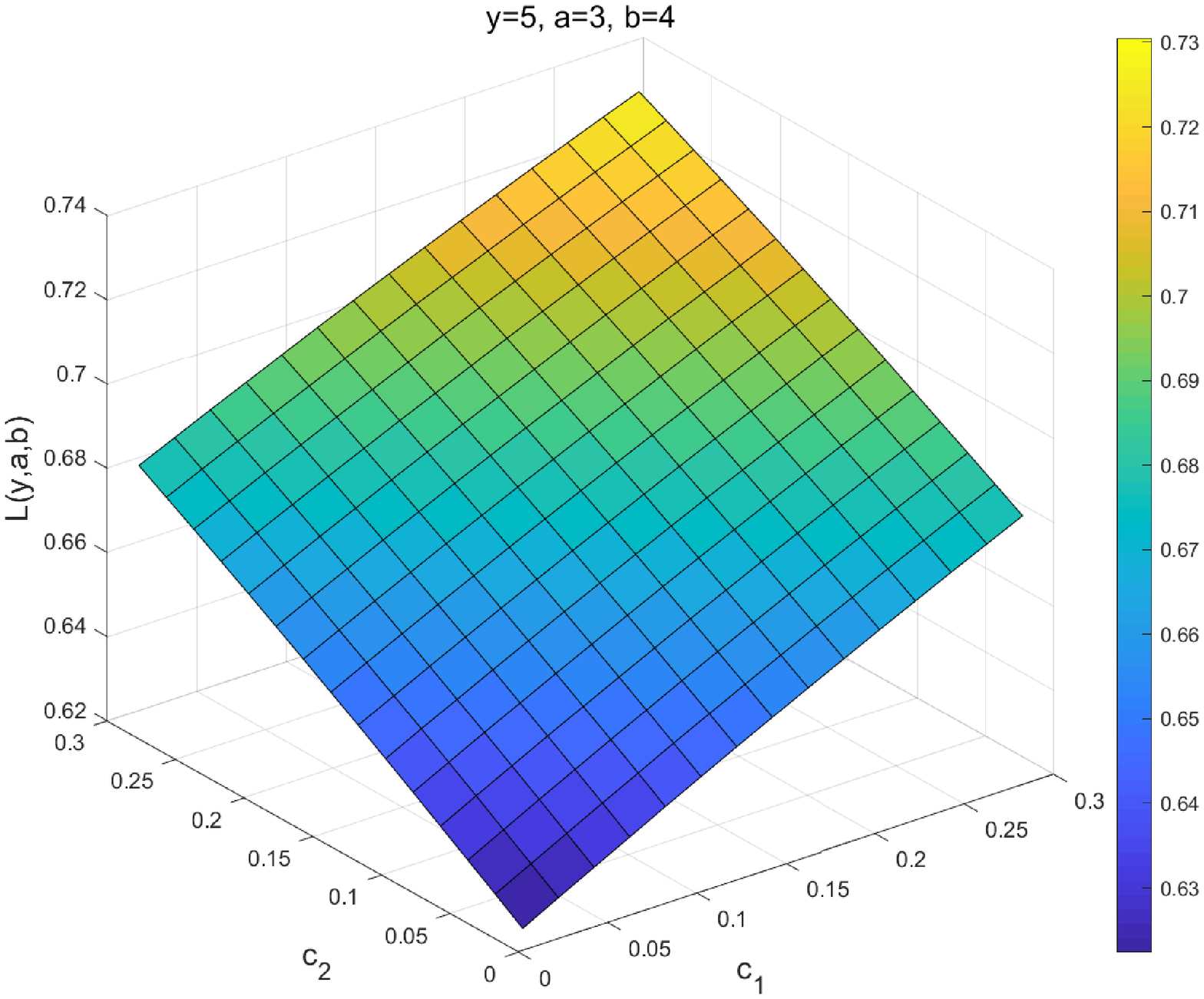}
}%
\caption{The Laplace transform of the ruin time $L(y;a,b)$ respectively as a function of
$\gamma$, $a$, $b$, $c_1$, $c_2$:(a) $y=8;a=3;b=4 $, (b) $y=8;b=5 $(solid) and $y=2;b=5 $(dashed), (c) $y=10;a=3$, (d) $y=5;a=3;b=4 $.
}
\label{fig.6}
\end{figure}

In this subsection, we focus on the Laplace transform of the ruin time. Also set the same basic parameter settings with previous subsection $\ref{ANA_div}$.
In Figure $\ref{fig.6}$, we characterize the effects of the main parameters $\gamma$, $a$, $b$, $c_1$, and $c_2$ on the Laplace transform of the ruin time $L(y;a,b)$.

First, we show in Figure $\ref{fig.6}$(a) that the curve of $L(y;a,b)$ as a function of the inter-dividend-decision times parameter $\gamma$ when $y=8$, $a=3$, and $b=4$, in which we denote by $L(8;3,4)$. We find that $L(8;3,4)$ is a monotonically increasing and convex function as $\gamma$. There are two factors that bring about this phenomenon. On one hand, smaller $\gamma$ means that more dividend decision times per unit time, which will lead to ruin earlier. On the other hand, when $\gamma$ is larger, dividend payments may become small, which will lead to the ruin time delayed.

Next, in Figure $\ref{fig.6}$(b-c), we plot function $L(y;a,b)$ with respect to barriers $a$ and $b$ under two different cases $y\ge b$ and $y<b$. From Figure $\ref{fig.6}$(b), we know that both $L(8;a,5)$ and $L(2;a,5)$ decrease as $a$ increases. We see also from Figure $\ref{fig.6}$(c) that $L(10;3,b)$ is a non-increasing function as $b$. We can explain these phenomena as follows. Larger $a$ or $b$ means that small potential dividends will be paid off, therefore ruin postpones. In the end, we plot $L(y;a,b)$ as a function of $(c_1,c_2)$ in Figure $\ref{fig.6}$(d). As can be seen from Figure $\ref{fig.6}$(d), one can observe that $L(5;3,4)$ is a decreasing function of both $c_1$ and $c_2$. When $c_1$ or $c_2$ is larger, more potential dividends will be paid off, thus the ruin time will be prolonged.

At the end of this subsection, it is worthwhile to mention an interesting fact. Comparing Figure $\ref{fig.1}$-$\ref{fig.4}$ and Figure $\ref{fig.6}(b-d)$ cross the same abscissa parameters, $V(y;a,b)$ appears to decrease while $L(y;a,b)$ increases.

\section{Conclusions}\label{Con}
In this paper, we studied dividend problems and ruin problems for spectrally negative L\'{e}vy risk model with ratcheting-periodic dividend strategy. Dividend payments can both be made discretely at the jump times of an independent Poisson process and continuously without falling.
The precise expressions of the expected NPV of dividends paid up to ruin and the Laplace transform of the ruin time are derived by using L\'{e}vy fluctuation theory and are written concisely in terms of scale functions.

Finally, we describe the two functions, the expected NPV of dividends and the Laplace transform of the ruin time, under Brownian risk model, a special spectral negative L\'{e}vy process. The optimal dividend value and the minimum Laplace transform value and the corresponding optimal barriers in the fixed settings are obtained. The two results are consistent as follows: $a^*=b$ and $b^*=\infty$. If we fixed the two barriers $a$, $b$, and other parameters, the optimal $c_1^*$ and $c_2^*$ should be $zero$. This indicates that if we do not consider the factor the investors are unwilling to see a decline in the dividend rate, the optimal mixed dividend strategy is pure periodic dividend strategy.

\vskip 3mm

%
%
\section*{References}
\vskip0.1in\parskip=0mm \baselineskip 15pt
\renewcommand{\baselinestretch}{1.15}
\footnotesize
\noindent
Albrecher, H., Azcue, P., Muler, N. 2020a. Optimal ratcheting of dividends in a Brownian risk model. arXiv:2012.10632.
\vskip 3mm

\noindent
Albrecher, H., Azcue, P., Muler, N. 2020b. Optimal ratcheting of dividends in insurance. {\it SIAM J. Control Optim.} 58,
1822--1845.
\vskip 3mm

\noindent
Albrecher, H., B\"{a}uerle, N., Bladt, M. 2018. Dividends: From refracting to ratcheting. {\it Insur. Math. Econ.} 83, 47--58.
\vskip 3mm

\noindent
Avanzi, B., Lau, H., Wong, B. 2020a. Optimal periodic dividend strategies for spectrally positive L\'{e}vy risk processes
with fixed transaction costs. {\it Insur. Math. Econ.} 93, 315--332.
\vskip 3mm

\noindent
Avanzi, B., Lau, H., Wong, B. 2021. On the optimality of joint periodic and extraordinary dividend strategies. {\it Eur. J.
Oper. Res.} 295, 1189--1210.
\vskip 3mm

\noindent
Avanzi, B., P\'{e}rez, J.L., Wong, B., Yamazakid, K. 2017. On optimal joint reflective and refractive dividend strategies in
spectrally positive L\'{e}vy models. {\it Insur. Math. Econ.} 72, 148--162.
\vskip 3mm

\noindent
Avanzi, B., Tu, V., Wong, B. 2016. On the interface between optimal periodic and continuous dividend strategies in the
presence of transaction costs. {\it ASTIN Bull.} 46, 708--745.
\vskip 3mm

\noindent
Avanzi, B., Tu, V., Wong, B. 2020b. Optimality of hybrid continuous and periodic barrier strategies in the dual model.
{\it Appl. Math. Optim.} 82, 105--133.
\vskip 3mm

\noindent
Avram, F., P\'{e}rez, J.L., Yamazaki, K. 2018. Spectrally negative L\'{e}vy processes with Parisian reflection below and classical
reflection above. {\it Stoch. Process. Their Appl.} 128, 255--290.
\vskip 3mm

\noindent
Chan, T., Kyprianou, A.E., Savov, M. 2011. Smoothness of scale functions for spectrally negative L\'{e}vy processes. {\it Probab.
Theory Relat. Field} 150, 691--708.
\vskip 3mm

\noindent
Dong, H., Yin, C., Dai, H. 2019. Spectrally negative L\'{e}vy risk model under Erlangized barrier strategy. {\it J. Comput.
Appl. Math.} 351, 101--116.
\vskip 3mm

\noindent
Dong, H., Zhou, X. 2019. On a spectrally negative L\'{e}vy risk process with periodic dividends and capital injections. {\it Stat.
Probab. Lett.} 155, 108589.
\vskip 3mm

\noindent
Kuznetsov, A., Kyprianou, A.E., Rivero, V. 2013. The theory of scale functions for spectrally negative L\'{e}vy processes.
Springer Berlin Heidelberg, Berlin, Heidelberg. pp. 97--186.
\vskip 3mm

\noindent
Kyprianou, A.E., Loeffen, R. 2010. Refracted L\'{e}vy processes. {\it Annales de l'Institut Henri Poincar\'{e}, Probabilit\'{e}s et
Statistiques} 46, 24--44.
\vskip 3mm

\noindent
Kyprianou, A.E., Loeffen, R., P\'{e}rez, J.L. 2012. Optimal control with absolutely continuous strategies for spectrally
negative L\'{e}vy processes. {\it J. Appl. Probab.} 49, 150--166.
\vskip 3mm

\noindent
Li, P., Meng, Q., Yuen, K.C., Zhou, M. 2021. Optimal dividend and risk control policies in the presence of a fixed
transaction cost. {\it J. Comput. Appl. Math.} 388, 113271.
\vskip 3mm

\noindent
Liu, Z., Chen, P., Hu, Y. 2020. On the dual risk model with diffusion under a mixed dividend strategy. {\it Appl. Math.
Comput.} 376, 125115.
\vskip 3mm

\noindent
Loeffen, R. 2008. On optimality of the barrier strategy in de Finetti's dividend problem for spectrally negative L\'{e}vy
processes. {\it Ann. Appl. Probab.} 18, 1669--1680.
\vskip 3mm

\noindent
Noba, K., P\'{e}rez, J.L., Yamazaki, K., Yano, K. 2018. On optimal periodic dividend strategies for L\'{e}vy risk processes.
{\it Insur. Math. Econ.} 80, 29--44.
\vskip 3mm

\noindent
P\'{e}rez, J.L., Yamazaki, K. 2018. Mixed periodic-classical barrier strategies for L\'{e}vy risk processes. {\it Risks} 6, 33.
\vskip 3mm

\noindent
Shen, Y., Yin, C., Yuen, K.C. 2013. Alternative approach to the optimality of the threshold strategy for spectrally
negative L\'{e}vy processes. {\it Acta Math. Appl. Sin.-Engl. Ser.} 29, 705--716.
\vskip 3mm

\noindent
Song, Z., Sun, F. 2021. The dual risk model under a mixed ratcheting and periodic dividend strategy. {\it Commun. Stat.
Theory Methods}, 1--15.
\vskip 3mm

\noindent
Yin, C., Shen, Y., Wen, Y. 2013. Exit problems for jump processes with applications to dividend problems. {\it J. Comput.
Appl. Math.} 245, 30--52.
\vskip 3mm

\noindent
Yin, C., Wen, Y. 2013. Optimal dividend problem with a terminal value for spectrally positive L\'{e}vy processes. {\it Insur.
Math. Econ.} 53, 769--773.
\vskip 3mm

\noindent
Yin, C., Wen, Y., Zhao, Y. 2014. On the optimal dividend problem for a spectrally positive L\'{e}vy process. {\it ASTIN Bull.
44, 635--651.
\vskip 3mm

\noindent
Yin, C., Yuen, K.C. 2011. Optimality of the threshold dividend strategy for the compound Poisson model. {\it Stat. Probab.
Lett.} 81, 1841¨C1846.
\vskip 3mm

\noindent
Zhang, A., Liu, Z. 2020. A L\'{e}vy risk model with ratcheting dividend strategy and historic high-related stopping. {\it Math.
Probl. Eng.} 2020, 6282869.
\vskip 3mm

\noindent
Zhang, Z., Han, X. 2017. The compound Poisson risk model under a mixed dividend strategy. {\it Appl. Math. Comput.} 315, 1--12.
\end{document}